\documentclass{elsart}
\usepackage{mathrsfs}
\usepackage{latexsym}
\usepackage{amssymb}
\usepackage{amsmath}
\usepackage{amsfonts}
\usepackage{indentfirst}
\usepackage{color}

\usepackage{epsfig,psfrag}

\usepackage{graphicx}
\usepackage{epstopdf}

\usepackage{subfigure}

\allowdisplaybreaks

\renewcommand\thesection{\arabic{section}}
\renewcommand{\thefootnote}{\arabic{footnote}}
\newtheorem{theorem}{\color{black}\indent Theorem}
\newtheorem{lemma}{\color{black}\indent Lemma}
\newtheorem{proposition}{\color{black}\indent Proposition}
\newtheorem{definition}{\color{black}\indent Definition}
\newtheorem{remark}{\color{black}\indent Remark}
\newtheorem{corollary}{\color{black}\indent Corollary}
\newtheorem{example}{\color{black}\indent Example}

\journal{}

\begin{document}

\begin{frontmatter}

\title{Differential Galoisian approach to Jacobi integrability of general analytic dynamical systems and its application}

\author{Kaiyin Huang$^{a}$},
\author{Shaoyun Shi$^{b,c,*}$},
\author{Shuangling Yang$^{b}$},
\thanks{ Corresponding author *.\\  telephone numbers: +86-431-85168423.
\\\textit{Email address}:
 huangky@scu.edu.cn(K.Huang), yangsl19@mails.jlu.edu.cn(S.Yang), shisy@jlu.edu.cn(S.Shi)   }

\address{$^a$ School of Mathematics, Sichuan University, Chengdu 610000,  China\\
$^b$ School of Mathematics, Jilin University, Changchun 130012,  China\\
$^c$ State Key Laboratory of Automotive Simulation and Control, Jilin University, Changchun 130012, P. R. China\\
}

\begin{abstract}
The Morales-Ramis theory provides an effective and powerful non-integrability criterion for complex
analytical Hamiltonian systems \emph{via} the differential Galoisian obstruction. In this paper we give a new
Morales-Ramis type theorem on the meromorphic Jacobi non-integrability of general analytic dynamical systems.
  The key point is to show the existence of Jacobian multiplier of a nonlinear system implies the existence of common Jacobian multiplier of Lie algebra associated with the identity component.
   In addition, we apply our results to the polynomial integrability of Karabut systems for stationary gravity waves in finite depth.



\end{abstract}
\end{frontmatter}

\section{Introduction}

The fundamental problem in the field of dynamical systems is to detect whether a given system is integrable or not.
 Roughly speaking, a system is called integrable if it has a number of invariant tensors(for instance, first integrals, symmetry fields or Jacobian multipliers) such that it can be solved by quadrature or in closed form.
The integrability structure of a system provides us its general solutions in an ``explicit'' way, and helps us
obtain the global dynamics, topological structure or final evolution of phase curves for the considered system \cite{gg,zx,in1,in2,in3}. On the contrary, non-integrability of a system pushes us  to expect that the system admits chaotic phenomena or
complex dynamical behaviors \cite{non1,non2,non3}.

It should be pointed out that there is no unique definition of integrability for dynamical systems.
In terms of Hamiltonian systems, the integrability is well defined in the Liouville sense, that is,
an $n$-degree of freedom
Hamiltonian system is integrable if and only if it has $n$ functionally independent first integrals in involution.
In terms of general non-Hamiltonian systems, there exist several different definitions of integrability as follows.

Let us recall some basic notions and facts concerning integrability of  general non-Hamiltonian systems.
Consider a general analytic system of differential equations
\begin{equation}
\dot x=F(x),~~~~x=(x_1,\cdots,x_n)\in \mathbb{C}^n,\label{ae}
\end{equation}
with smooth right-hand sides $F=(F_1,\cdots,F_n)$. Denote by $\mathcal{L}_F$ the Lie derivative
associated with system (\ref{ae}).
A differential tensor field $T(x)$ is called an invariant tensor of system (\ref{ae}) if $\mathcal{L}_F(T)=0$.
Three typical invariant tensors are first integrals, symmetry field and invariant $n$-form. More precisely, a scalar function $\Phi(x)$ is a first integral of (\ref{ae}) if and only if $\mathcal{L}_F(\Phi):=\langle \partial_x \Phi,F \rangle=0$, where $\partial_x \Phi$ is the gradient of $\Phi$ and $\langle,\rangle$ denotes the
inner product in $\mathbb{C}^n$. An $n$-dimensional vector function $V(x)$ is a symmetry field of (\ref{ae}) if and only if $\mathcal{L}_F(V):=[V,F]=0$. Clearly, system (\ref{ae}) always has a \emph{trivial} symmetry field, that is, itself $V=F$.
 An $n$-form $\Omega=J(x)dx_1\wedge \cdots\wedge dx_n$ is an invariant $n$-form of (\ref{ae}) if and only if $\mathcal{L}_F(\Omega):=div(JF)dx_1\wedge \cdots\wedge dx_n=0$, where the scalar function $J(x)$
is called a Jacobian multiplier of (\ref{ae}).

System (\ref{ae}) is completely integrable if it admits $(n-1)$ functionally independent first integrals $\Phi_1,\cdots,\Phi_{n-1}$. Clearly, the orbits of a completely integrable system are contained in the  curves
$$
\mathcal{S}_{c_1\cdots c_{n-1}}=\{x|\Phi_1(x)=c_1,\cdots,\Phi_{n-1}=c_{n-1}\}.
$$
Let us mention that system with $(n-1)$ $C^r$ first integrals is orbitally
equivalent to a linear differential system in a full Lebesgue subset \cite{js2016}.

Another definition of integrability is due to the classical work of S. Lie, who proved that system admitting $n$ linearly independent and commuting symmetries $V_1=F,V_2,\cdots,V_n$ is integrable by quadrature \cite{jo2}.
Correspondingly, system (\ref{ae}) is called integrable in the Lie sense if it has $n$ linearly independent symmetries $V_1=F,V_2,\cdots,V_n$ such that $[V_i,V_j]=0$ for any $1\leq i,j\leq n$.

Doing a careful analysis of the concept of the Liouville integrability of Hamiltonian systems,
Bogoyavlenskij \cite{B} proposed a new definition of integrability for non-Hamiltonian systems, which is regarded as a generalization of the Liouville integrability. More precisely, system (\ref{ae}) is integrable
in the Bogoyavlenskij sense if for some $k\in\{0,1,\cdots,n-1\}$, it has $k$ functionally independent
first integrals $\Phi_1,\cdots,\Phi_k$ and $n-k$ linearly independent vector fields $V_1=F,\cdots,V_{n-k}$ such that
$$
[V_i,V_j]=0,~~\text{and}~~\langle\partial_x\Phi_l,V_j\rangle=0,~~\text{for}~~1\leq l\leq k,~~1\leq i,j\leq n-k.
$$
Similar to the case of Hamiltonian systems, if system (\ref{ae}) is integrable in the Bogoyavlenskij
sense, then the invariant sets associated with first integrals
are generically diffeomorphic to tori, cylinders, or planes inside the phase space \cite{B}.
Obviously, complete integrability and Lie integrability are special cases of  Bogoyavlenskij integrability, corresponding to $k=0$ or $n-1$ respectively.

The next definition of integrability is due to Jacobi \cite{jo0}, and is widely applied for nonholonomic mechanical systems \cite{jo1,jo2}. System (\ref{ae}) is called integrable in the Jacobi sense if it has $n-2$ functionally independent
first integrals $\Phi_1,\cdots,\Phi_{n-2}$ and a Jacobian multiplier $J(x)$. The integrability in terms of the
existence of $n-1$ Jacobian multipliers was studied in \cite{hu2014,zhang2014}.

In 2013, Kozlov combined the above definitions and proposed the Euler-Jacobi-Lie integrability: if system (\ref{ae}) has $k$ functionally independent first integrals $\Phi_1,\cdots,\Phi_k$, $(n-k-1)$ independent symmetry fields $V_1=F,V_2,\cdots,V_{n-k-1}$ generating a nilpotent Lie algebra of vector field, and an invariant volume $n$-form $\Omega=J(x)dx_1\wedge \cdots\wedge dx_n$ such that
$$
\mathcal{L}_{V_i}(\Phi_j)=0,~~\mathcal{L}_{V_i}(\Omega)=0,~~1\leq i\leq n-k-1,~~1\leq j\leq k,
$$
then system (\ref{ae}) can be integrated by
quadrature \cite{kozlov2013}. To summarize,
all above accepted definitions of the integrability follow the same philosophy: the existence of some invariant tensors, including first integrals, symmetry fields and Jacobian
multipliers, whose total number of invariant tensors is equal to the dimension $n$ of
the dynamical system.
\begin{table}[htp]
\centering
\caption{Summary of integrability for differential systems in terms of invariant tensors.}
\begin{tabular}{|c|c|c|c|}
\hline
Integrability &First integrals& Symmetry fields& Jacobian multipliers   \\
\hline
Complete integrability  &$n-1$ &$1$ &0\\

Lie integrability  & 0&$n$ &0\\

Bogoyavlenskij integrability  &$k$&$n-k$ & 0\\

Jacobi integrability &$n-2$ & 1&1\\

Integrability in \cite{hu2014,zhang2014} &0 & 1&$n-1$\\

Euler-Jacobi-Lie integrability  & $k$&$n-1-k$& 1\\

\hline
\end{tabular}
\end{table}

There are few effective methods to decide whether a system is integrable or not.
 At the end of the 20th century, a significant development, called the Morales-Ramis theory, was made by  Morales, Ramis, Sim\'{o} \cite{ramis,morales2001,morales1994picard}, and Baider,
Churchill, Rod, Singer \cite{baider1996,churchill}, who used the properties of differential Galois group of variational equations to give strong and effective necessary
conditions for the  Liouville integrability of Hamiltonian systems. This theory has been  applied successfully for a number of nonlinear physical models such as N-body problem \cite{ga1,ga11}, Hill's problem \cite{ga2}, problems with homogeneous  potentials \cite{ga3,ga4}, geodesic motions \cite{ga5} and other physical problems \cite{ga6,ga7,ga71}.
 In the recent decades, inspired by the works of Morales \emph{et al.},
the differential Galoisian approach
has been attempted to the integrability analysis of non-Hamiltonian systems. In 2010, Ayoul and Zung \cite{zung} used
the cotangent lifting trick to naturally extend the Morales-Ramis theory into non-Hamiltonian systems, and gave
  necessary conditions for meromorphic integrability  in the Bogoyavlenskij sense. The necessary conditions for the complete integrability are proposed in \cite{li2012,li2008}, which are expressed in terms of finiteness of the differential Galois group. Based on Malgrange pseudogroup and Artin approximation,  Casale \cite{jacobi1} showed that if
  a rational vector field is rationally integrable in the Jacobi sense on an
algebraic variety then identity components of Galois groups of variational
equations are solvable and their first derived Lie algebras are abelian, see \cite{jacobi2,jacobi3}  for more details
on Jacobi non-integrability.

   Our main results of this work are summarized as follows.
\begin{itemize}
\item[(i)]
Instead of the tools of Malgrange pseudogroup and Artin approximation used in Casale's proof, we provide an
 elementary proof of necessary conditions for the Jacobi integrability in the category of meromorphic functions. Indeed, we show that if system (\ref{ae}) has $0\leq k\leq n-2$ meromorphic first integrals $\Phi_1,\cdots,\Phi_k$ and $n-1-k$ meromorphic Jacobian multipliers $J_1,\cdots,J_{n-k-1}$ such that
$$
\Phi_1,\cdots,\Phi_k,\frac{J_2}{J_1},\cdots,\frac{J_{n-1-k}}{J_1},
$$
are functionally independent, then the identity
component of the differential Galois group of the  normal variational equations along a particular solution is abelian (\emph{Theorem \ref{thbb}}).  Our proof strategy can be used to investigate necessary conditions for other integrability via differential Galoisian methods.

\item[(ii)]
Combining with Theorem \ref{thbb}, we study the polynomial integrability of  Karabut systems for stationary gravity waves in finite depth.
 Witting \cite{wi} proposed  a new formal series solution to the water waves problem when he studied the solitary wave in a fluid of finite depth.
Then Karabut \cite{ka1,ka2,ka3} showed
the problem of exact summation of Witting's series can be reduced to solve or integrate some homogeneous ordinary differential equations, called Karabut systems. We show that the 3-dimensional
Karabut system is integrable and admits infinitely many Hamilton-Poisson
realizations and a Lax formulation (\emph{Proposition \ref{pro1},\ref{pro2},\ref{pro3}}), and also show that the 5-dimensional Karabut system has two and only two functionally independent
polynomial first integrals (\emph{Theorem \ref{mth}}) which answers the question by Karabut \cite{ka3} and improves the result in \cite{k1} from the point
of view of partial integrability.

\end{itemize}

The outline of this paper is as follows. In section 2, we introduce some preliminary
notions and results on the differential Galois theory and Jacobian multipliers. In section 3, we
formulate a general theorem which provides necessary conditions for the existence of $k$ first integrals and $n-1-k$ Jacobian multipliers satisfying compatibility conditions. In section 4, we apply our results to study integrability of  3-dimensional and 5-dimensional Karabut systems.

\section{Preliminaries results}

\subsection{Differential Galois theory}
Differential Galois theory is a  generalization of the classical Galois theory from  polynomial equations to linear differential equations,
 which is also called Picard-Vessiot theory.  Here we only
review some parts of differential Galois theory, for
more details see \cite{ramis,van} and the references therein.

Consider the linear homogeneous differential equations on a
differential field $(\mathbb{K}, \partial)$
\begin{equation}\label{q2}
Y'=AY, \ \ \ A\in Mat(\mathbb{K}, n),
\end{equation}
where  $Mat(\mathbb{K},n)$ denotes the ring of $n\times n $ matrices with entries in $\mathbb{K}$.
Recall that a differential field is a pair $(\mathbb{K}, \partial)$ consisting of a field
$\mathbb{K}$ and a derivative $\partial$, where the derivative $\partial$ is an additive mapping
$\partial: \mathbb{K}\rightarrow
 \mathbb{K}$ satisfying
\begin{equation*}\label{q1}
\partial(ab)=\partial(a)b+a\partial(b), \ \ \ a, b\in \mathbb{K}.
\end{equation*}
We also write $a'$ instead of
$\partial(a)$. The set of elements of $\mathbb{K}$ for which $\partial$ vanishes is called the field of constants of $\mathbb{K}$,
denoted by ${\rm Const}(\mathbb{K}):=\{a\in\mathbb{K}|a'=0\}$.
In practical applications,  $\mathbb{K}$ is the field of meromorphic functions on a Riemann surface endowed with a meromorphic vector field and the field of constants becomes the field of complex numbers $\mathbb{C}$.

 A differential field extension,
denoted by $\mathbb{L} / \mathbb{K}$, is a field extension such that $\mathbb{L}$ is also a differential field and the derivations on $\mathbb{L}$ and $\mathbb{K}$ coincide on $\mathbb{K}$, i.e., $\partial_\mathbb{L}|_{\mathbb{K}}=\partial_\mathbb{K}$, where $\partial_\mathbb{L}, \partial_\mathbb{K}$
are the derivation of $\mathbb{L}, \mathbb{K}$ respectively. Now we introduce two important differential field extensions. The first one is the Liouvillian extension.
\begin{definition}\label{d3}
The differential field extension $\mathbb{L}/ \mathbb{K}$ is called a Liouvillian extension if  ${\rm Const}(\mathbb{L} )={\rm Const}(\mathbb{K} )$ and there exists a tower of extensions
$$
\mathbb{K}=\mathbb{L}_0 \subset \mathbb{L}_1 \subset \cdots \subset \mathbb{L}_m=\mathbb{L},
$$
such that for $i=1, \cdots, m,~ \mathbb{L}_i=\mathbb{L}_{i-1}(a_i)$, and one of the
following cases holds:

$(1)$ $a_i'\in \mathbb{L}_{i-1}$, in this case  $a_i$ is called an integral element
of $L_{i-1}$;

$(2)$ $a_i\neq0$ and $a_i'/a_i\in \mathbb{L}_{i-1}$, in this case $a_i$ is
called an exponential integral element of $\mathbb{L}_{i-1}$;

$(3)$ $a_i$ is algebraic over $\mathbb{L}_{i-1}$.
\end{definition} 
\begin{remark}
Roughly speaking, the fact that $\mathbb{L}/ \mathbb{K}$ is a Liouvillian extension means that each element of  $\mathbb{L}$ can be  built up from $\mathbb{K}$ by algebraic operations and taking exponentials or
indefinite integrals.
\end{remark}

The second one is  the Picard-Vessiot(P-V) extension which is associated with a linear system of differential equations (\ref{q2}).

\begin{definition}\label{d3}
The differential field extension  $\mathbb{L}/\mathbb{K}$ is a Picard-Vessiot(P-V) extension for the linear system $(\ref{q2})$ if and only if it satisfies the following three conditions:

$(1)$ ${\rm Const}(\mathbb{L} )={\rm Const}(\mathbb{K} )$;

$(2)$ There exists a fundamental matrix $\Phi\in GL(\mathbb{L}, n)$ for the linear homogeneous differential equations
$(\ref{q2})$;

$(3)$ $\mathbb{L}$ is generated over $\mathbb{K}$ as a differential field by the
entries of the fundamental matrix $\Phi$.
\end{definition} 

\begin{remark}
Roughly speaking, the P-V extension is the smallest differential extension such that it contains $n$ linearly independent solutions of (\ref{q2}) and no new constants are added. In addition,  it is known that if the constant subfield of differential field $(\mathbb{K},
\partial)$ is characteristic zero,  for instance
${\rm Const}(\mathbb{K})=\mathbb{C}$, then (\ref{q2}) can
admit a Picard-Vessiot extension which is unique up to isomorphism \cite{van}.

\end{remark}

Fixed a P-V extension $\mathbb{L}/\mathbb{K}$ and the fundamental matrix $\Phi$,
all differential $\mathbb{K}$-automorphisms ($\sigma:\mathbb{L}\rightarrow \mathbb{L} ,
 \sigma(a')=(\sigma(a))', \forall a \in \mathbb{L} $ and $\sigma(a)=a, \forall a\in
 \mathbb{K}$) of $\mathbb{L}$ is called the differential Galois group of (\ref{q2}) and denoted by $Gal(\mathbb{L}/ \mathbb{K})$.
 Let $\Phi(t)$ be a fundamental-solution
matrix of (\ref{q2}). Note that for any $\sigma \in Gal(\mathbb{L}/\mathbb{K})$, we have $\sigma(\Phi)$ is also a fundamental
matrix of (\ref{q2}). Therefore $\sigma(\Phi)=\Phi M_\sigma$ with
$M_\sigma\in GL(\mathbb{C}, n)$, which gives a faithful representation of the group of $\mathbb{K}$-automorphisms of $\mathbb{L}$ on the general linear group as
$$
 \rho:  Gal(\mathbb{L}/\mathbb{K})\rightarrow GL(\mathbb{C}, n),~~\sigma\rightarrow M_\sigma.
$$
Hence, the differential Galois group $Gal(\mathbb{L}/\mathbb{K})$ can be regarded as a subgroup of
$GL(\mathbb{C},n)$.  In what follows, we will be dealing either with the differential Galois group $G:=Gal(\mathbb{L}/\mathbb{K})$
or its matrix group $\rho(G)$. Moreover, $G$ is a linear algebraic group \cite{van}, and is a union of a finite number of disjoint connected
components.  Then the differential Galois group $G$ has a unique maximal connected subgroup $G^0$  containing the identity element of the group, called the identity component of $G$. We say that a group $G$ is solvable if and only if there exists a chain of normal subgroups
$$
e=G_{0} \triangleleft G_{1} \triangleleft \ldots \triangleleft G_{n}=G,
$$
such that  the quotient $G_{i+1}/G_i$ is abelian for $i=0,\cdots,n-1$.

The following fundamental theorem concerns the deep relation between the solvability of the linear differential equations and  the corresponding differential
Galois group.

\begin{theorem}\label{t3}
Let  $\mathbb{L}/ \mathbb{K}$ be a P-V extension for $(\ref{q2})$. Then the linear system $(\ref{q2})$ is solvable
by quadrature, that is, $\mathbb{L}/ \mathbb{K}$ is a Liouvillian extension if and only if the identity component $Gal(\mathbb{L}/\mathbb{K})^0$ of differential Galois
group $Gal(\mathbb{L}/\mathbb{K})$ is solvable. In particular, if  the
identity component $Gal(\mathbb{L}/\mathbb{K})^0$ is abelian, (\ref{q2})
is solvable by quadrature.
\end{theorem}

In general, for an $n$-dimensional linear differential system, it is difficult to compute or analyze
the properties of the corresponding differential Galois group. When the (normal)    variational equations can
be reduced into a second-order linear differential equation with rational coefficients, the so-called
Kovacic's algorithm \cite{k} can help us calculate  effectively the solvability of the differential Galois group.

\begin{lemma}\cite{k}\label{kth}
The differential Galois group $G$ of
\begin{equation}\label{nnve}
\frac{d^2\chi}{dt^2}=r(t)\chi,~~~r(x)\in \mathbb{C}(t),
\end{equation}
with $\mathbb{C}(t)$ being the  rational function field on $\mathbb{C}$
can be classified with the following four cases.

{\bf Case 1.} $G$ is conjugate to a subgroup of triangular group
\begin{equation}\label{tri}
\mathbb{C}^{*} \ltimes \mathbb{C}^{+}=\left\{\left(\begin{array}{cc}
c & 0 \\
b & c^{-1}
\end{array}\right) \mid b \in \mathbb{C}, c \in \mathbb{C}^*\right\}.
\end{equation}
Then $(\ref{nnve})$ has a solution of the form $e^{\int \omega}$
with $\omega\in\mathbb{C}(t) $.

{\bf Case 2.}  $G$ is not of case 1, but is conjugate to a subgroup of the infinite dihedral group
\begin{equation*}
D = \left\{ \left( {\begin{array}{*{20}{c}}
c&0\\
0&{{c^{ - 1}}}
\end{array}} \right)\left| c\in \mathbb{C}^* \right.\right\}  \cup \left\{ \left( {\begin{array}{*{20}{c}}
0&c\\
{ - {c^{ - 1}}}&0
\end{array}} \right)\left| c\in \mathbb{C}^* \right.\right\}.
\end{equation*}
Then $(\ref{nnve})$ has a solution of the form $e^{\int \omega}$
with $\omega$  algebraic over $\mathbb{C}(t)$ of degree two.

{\bf Case 3.} $G$ is not of case 1 and case 2, but is a finite group. Then all solutions of $(\ref{nnve})$ are algebraic over~$\mathbb{C}(t)$.

{\bf Case 4.} $G=SL(2,\mathbb{C})$, where $SL(2,\mathbb{C})$ is  the group of $2\times2$ matrices with elements in $\mathbb{C}$ and determinant one. Then $(\ref{nnve})$ is not integrable in Liouville sense.
\end{lemma}

{\bf Proof:}
See Lemma 1.4 in \cite{k}.\qed

\begin{remark}
Kovacic \cite{k} presented a complete algorithm to analyze which cases the differential Galois group $G$ of system (\ref{nnve}) falls into.
Due to Lemma \ref{kth}, there exist four cases in Kovacic's algorithm.
We can obtain the Liouvillian solutions of (\ref{nnve}) in
the cases 1,2,3, but for the case 4 this equation has no Liouvillian solutions. In addition,
Kovacic's algorithm gives us only one solution $\xi_1$  in the case 1,2,3, and the second independent solution
$\xi_2$ can be obtained by
$$
\xi_2=\xi_1\int\frac{dt}{\xi_1^2}.
$$
\end{remark}
\begin{remark}
For a general second-order linear differential equation
\begin{align}\label{nnnve}
 \frac{d^2\xi}{dt^2}+a(t) \frac{d\xi}{dt}+b(t)\xi=0,~~~a(t),b(t)\in \mathbb{C}(t),
\end{align}
we can make a well-known change of
variable
$$
\xi=\chi\exp\bigg(-\frac{1}{2}\int_{{t_0}}^t {a(s)ds} \bigg),
$$
and get the reduced form of (\ref{nnnve})
\begin{equation}\label{nnve3}
\frac{d^2\chi}{dt^2}=r(t)\chi,~~~r(t)=\frac{a^2}{4}+\frac{1}{2}\frac{da}{dt}-b.
\end{equation}
It should be pointed out that the identity component of (\ref{nnnve}) is solvable if and only if that of (\ref{nnve3}) is solvable,
since the above transformation does not affect the Liouvillian solvability
of (\ref{nnnve}).
\end{remark}

The next lemma is due to the work of  Singer and
Ulmer \cite{singer1993}, which provides a more precise characterization of case 1 of Lemma \ref{kth}.
\begin{lemma}\label{sth}
Assume that the differential Galois group  $G$ of (\ref{nnve}) is conjugate to a subgroup of triangular group. Then

{\bf Case 1.1} $G$  is conjugate to a subgroup of the diagonal group
$$
\left\{\left(\begin{array}{cc}
c & 0 \\
0 & c^{-1}
\end{array}\right) \mid c\in \mathbb{C}^*\right\}.
$$
In this case system (\ref{nnve}) has two independent solutions $\xi_1$ and $\xi_2$ such that $\xi_i'/\xi_i\in\mathbb{C}(t)$,~$i=1,2$.

{\bf Case 1.2} $G$  is conjugate to the group
$$
\left\{\left(\begin{array}{cc}
c & 0 \\
b & c^{-1}
\end{array}\right) \mid b \in \mathbb{C}, c^{m}=1 \text { with an } m \in \mathbb{N}\right\},
$$
which is a proper subgroup of the group $\mathbb{C}^{*} \ltimes \mathbb{C}^{+}$.
 In this case system (\ref{nnve}) has only one solution $\xi$ (up to constant multiples) such that $\xi'/\xi\in\mathbb{C}(t)$, and $m$ is the smallest positive integer such that $\xi^m \in\mathbb{C}(t)$.

{\bf Case 1.3} $G$  is conjugate to the group $\mathbb{C}^{*} \ltimes \mathbb{C}^{+}$.
In this case system (\ref{nnve}) has only one solution $\xi$ (up to constant multiples) such that $\xi'/\xi\in\mathbb{C}(t)$, and  $\xi^m \notin\mathbb{C}(t)$ for any positive integer $m$.
\end{lemma}
{\bf Proof:}
See Proposition 4.2 in \cite{singer1993}.\qed

\begin{remark}
Based on the above results,  let us point out that:

(1) the identity component $G^0$ associated with (\ref{nnve})
is not solvable if and only if $G$ belongs to the case 4;

(2) the identity component $G^0$ associated with (\ref{nnve})
is not abelian  if and only if $G$ belongs to either the case 1.3 or the case 4.
\end{remark}

Acosta-Hum\'{a}nez and Bl\'{a}zquez-Sanz \cite{poly} gave a complete classification of differential Galois group of a second-order differential equation with polynomial coefficients.

\begin{lemma}\label{poly}\cite{poly}
Let $Q(t)\in \mathbb{C}[t]/\mathbb{C}$ be a  polynomial of degree $k>0$.
The differential Galois group $G$ of
\begin{equation*}
\frac{d^2\chi}{dt^2}=Q(t)\chi,~~~Q(t)\in \mathbb{C}[t]/\mathbb{C},
\end{equation*}
is either $SL(2,\mathbb{C})$ or $\mathbb{C}^{*} \ltimes \mathbb{C}^{+}$,  in particular is non-abelian.
\end{lemma}

In applications, the next two results  can help us reduce the dimension of linear differential equations.
\begin{lemma}\cite{jordan}\label{sub}
Consider the differential Galois group $G$ of the linear differential system
\begin{equation}\label{esys}
\frac{d}{dt}\left( {\begin{array}{*{20}{c}}
{X_1}\\
{X_2}
\end{array}} \right) = \left( {\begin{array}{*{20}{c}}
{A}&{0}\\
{B}&{C}
\end{array}} \right)\left( {\begin{array}{*{20}{c}}
{X_1}\\
{X_2}
\end{array}} \right),~~~~~~~A\in Mat(\mathbb{K}, m),~~C\in Mat(\mathbb{K}, l).
\end{equation}
The following statements hold.

(i) If the identity component of  $G$ is solvable(abelian), then the identity component of the differential Galois group of the subsystem
\begin{align}\label{subsys}
\frac{d}{dt}X_1=AX_1
\end{align}
is also solvable(abelian).

(ii) Suppose $B=0$. Then system (\ref{esys}) is the direct sum of two subsystems (\ref{subsys}) and
\begin{align}\label{subsys2}
\frac{d}{dt}X_2=CX_2.
\end{align} The identity component of  $G$ is solvable(abelian) if and only if  the identity components of the differential Galois group of both  subsystems (\ref{subsys}) and (\ref{subsys2})
are  solvable(abelian).
\end{lemma}

The next result is a variant of the well-known Ziglin lemma \cite{ziglin} for which the meromorphic functions are replaced by
$C^2$ functions, but it is very useful to reduce the dimension of linear differential equations.
\begin{lemma}\label{ss}
Assume $\psi(t)$ is a particular solution of the nonlinear system (\ref{ae}).
If a $C^2$ function $\Phi(x)$ is a first integral of nonlinear system (\ref{ae}), then the function $G(t,\xi)=\langle\nabla \Phi(\psi(t)), \xi\rangle$
is either a constant or a time-dependent linear first integral of the  variational system
\begin{equation}\label{ave00}
\frac{{ d\xi}}{{ d}t}=A(t)\xi,~~A(t)={\frac{\partial F}{\partial x}}\bigg|_{x=\psi(t)}.
\end{equation}
\end{lemma}
{\bf Proof:}
By definition, we have $d\Phi(x(t))/dt=0$ for all solutions $x(t)$ of (\ref{ae}), that is,
\begin{align}\label{ada}
\langle\partial_x \Phi(x), F(x)\rangle=\sum_j{F_j\frac{\partial \Phi}{\partial x_j}}\equiv 0.
\end{align}
Taking the derivative of (\ref{ada}) with respect to $x_i$, we have
\begin{equation*}
\sum_j{\left(F_j\frac{\partial^2 \Phi}{\partial x_j\partial x_i}+\frac{\partial \Phi}{\partial x_j}\frac{\partial F_j}{\partial x_i}\right)}\equiv 0,~~~~i=1,\cdots,n.
\end{equation*}
Then,
\begin{align*}
\frac{dG}{dt}&=\frac{\partial G}{\partial t}+\sum_j\frac{\partial G}{\partial \xi_j}\frac{d\xi_j}{dt}\\
&=\sum_i\left(F_j\frac{\partial^2 \Phi}{\partial x_j\partial x_i}    \right)\xi_i+\sum_j\left(\frac{\partial \Phi}{\partial x_j}\sum_i(\frac{\partial F_j}{\partial x_i}\xi_i)\right)\\
&=\sum_i\left(\sum_j{\left(F_j\frac{\partial^2 \Phi}{\partial x_j\partial x_i}+\frac{\partial \Phi}{\partial x_j}\frac{\partial F_j}{\partial x_i}\right)}\right)\xi_i\\
&=0,
\end{align*}
which completes the proof.\qed

\subsection{(Generalized) Jacobian multipliers and the characterization}

The first result is due to the classical work of Poincar\'{e}, which gives an equivalent characterization of Jacobi multiplies. For a proof, please see Proposition 2.2 in \cite{zx}.

\begin{lemma}\label{jaco1}
Let $J(x)$ be a non-zero continuously differentiable function. Then the following statements are equivalent.

(i) $J(x)$ is a  Jacobian multiplier of system (\ref{ae}), that is,
$$
div(JF)=J(x)\nabla_x\cdot F(x)+\langle \partial_xJ,F\rangle\equiv0.
$$

(ii) For any  flow $\phi_t(t_0,x_0)$  of (\ref{ae}) satisfying $\phi_{t_0}(t_0,x_0)=x_0$, we have
$$
J(x_0)=J(\phi_t)\det \partial_{x_0}(\phi_t),~~\forall~t\geq t_0.
$$

(iii) For any bounded region $D_0$, the integral
$$
V(t)=\int_{D_t}J(y)dy.
$$
is independent with the time $t$, where $D_t=\{y|y=\phi_t(t_0,x_0),~x_0\in D_0\}$ is the evolution of $D_0$ under the flow  $\phi_t$.
\end{lemma}

Next we discuss the relationship between  Jacobian multipliers of equivalent differential systems.
\begin{lemma}\label{jaco2}
Assume $x=G(y)$ is a continuously differentiable and invertible transformation. If $J(x)$ is a Jacobian multiplier
of system (\ref{ae}), then
$$
\hat J(y):=J(G(y))\det \partial_yG(y)
$$
 is a Jacobian multiplier of system
$$
\dot y=(\partial_y G)^{-1}F(G(y)).
$$
\end{lemma}

Lemma \ref{jaco2} can be easily proved by using the definition of Jacobian multipliers or using  the invariant of integral $V(t)$, for a detailed proof we refer to Proposition 2.3 in \cite{zx}.

Now we introduce the notation of \emph{time-dependent} Jacobian multipliers for non-autonomous differential systems.
\begin{definition}
A non-zero  continuously differentiable  function $J(t,x)$ is called a Jacobian multiplier of non-autonomous differential systems
\begin{align}\label{nae}
\dot x=F(t,x),~~~x\in\mathbb{C}^n,
\end{align}
if
$$\partial_t J(t,x)+J(t,x)\nabla_x\cdot F(t,x)+\langle \partial_xJ,F\rangle\equiv0. $$
\end{definition}

If we set $\omega=t$ and rewrite (\ref{nae}) as an autonomous differential system
$$
\dot x=F(\omega,x),~~\dot \omega=1,
$$
then the above definition coincides with the usual Jacobian multiplier of autonomous differential systems. Furthermore,
as two corollaries of Lemmas \ref{jaco1}-\ref{jaco2}, we immediately obtain the following results.

\begin{lemma}\label{gjaco1}
Let $J(t,x)$ be a non-zero continuously differentiable function. Then $J(t,x)$ is a  Jacobian multiplier of system (\ref{nae}) if
and only if for  any  flow $\phi_t(t_0,x_0)$  of (\ref{nae}) satisfying $\phi_{t_0}(t_0,x_0)=x_0$, we have
$$
J(t_0,x_0)=J(t,\phi_t)\det \partial_{x_0}(\phi_{t}),~~\forall~t\geq t_0.
$$
\end{lemma}

\begin{lemma}\label{gjaco2}
Assume $x=G(t,y)$ is a continuously differentiable and invertible transformation for any fixed $t$. If $J(t,x)$ is a Jacobian multiplier
of system (\ref{nae}), then $\hat J(t,y):=J(t,G(t,y))\det \partial_yG(t,y)$ is a Jacobian multiplier of system
$$
\dot y=\hat F(t,y)=(\partial_y G)^{-1}\big(F(t,G(t,y))-\partial_tG(t,y)\big).
$$
\end{lemma}

\section{Necessary conditions for Jacobi integrability of analytic differential systems}

Let $\psi(t)$ be a non-equilibrium analytic solution of system (\ref{ae}).
Linearization of (\ref{ae}) around $\psi(t)$ yields the \emph{variational equations} of the following form:
\begin{equation}\label{ave}
\frac{{ d\xi}}{{ d}t}=A(t)\xi,~~A(t)={\frac{\partial F}{\partial x}}\bigg|_{x=\psi(t)}\in Mat(\mathbb{K},n),
\end{equation}
where  elements of field $\mathbb{K}$  are  meromorphic  functions on the phase curve $\Gamma = \{\psi(t)\}$.
Note that  $\dot \psi(t)$ is a
nontrivial solution of system (\ref{ave}) due to $\psi(t)$ being a non-equilibrium solution of system (\ref{ae}).
Using this fact, we can reduce the dimension of (\ref{ave}) by one.

Indeed, making a change of  variables
$$
\xi=P(t)\eta,
$$
where $P(t)=(P^1, \cdots, P^{n-1}, \dot{\psi}(t))$ is a non-singular
matrix with its components in $\mathbb{K}$, then (\ref{ave}) becomes the following equivalent form
\begin{equation}\label{e7}
\frac{{ d\eta}}{{ d}t}=B(t,\eta):=P(t)^{-1}(A(t)P(t)-\dot{P}(t))\eta=\begin{pmatrix}C(t) &\theta \\
\alpha(t)^T & 0\end{pmatrix}\eta,
\end{equation}
where $\theta$ denotes the $n-1$ dimensional zero vector and $\eta=(\zeta,\eta_1)^T$. Therefore, we obtain a subsystem of (\ref{e7})
\begin{align}\label{anve}
\frac{{ d\zeta}}{{ d}t}=C(t)\zeta,~~C(t)\in Mat(\mathbb{K},n-1),
\end{align}
which is  the so-called \emph{normal
variational equations} of (\ref{ae}) along $\Gamma$.

Equations $(\ref{ave})$, $(\ref{e7})$ and $(\ref{anve})$ are linear differential equations. Then, we can associate
the differential Galois theory with them.
Let us mention that the choice of $P(t)$ is not unique, and both (\ref{ave}) and (\ref{e7})
have the same
differential Galois group of equations since they have the same
P-V extension. Denote by $G$ the differential Galois group  of the normal  variational equations (\ref{anve}).
Recall that $G$ is a linear algebraic
group, thus, in particular, it is a Lie group, and one can
consider its Lie algebra which reflects only the properties
of the identity component $G^0$ of the group.
We denote by $\mathcal{G}\subset gl(n, \mathbb{C})$
the Lie algebra of $G$.
 Then arbitrary element $Y\in
\mathcal{G}$ can be viewed as a linear vector field: $x\rightarrow
Y(x):=Y\cdot x$, for $x\in \mathbb{C}^n$, and
$e^{Y t}\in G$ for all $t\in \mathbb{C}$.

The next result goes back to Ziglin \cite{ziglin} and plays a critical and fundamental result in the non-integrability approach, see Chapter 4 in the book \cite{ramis} for a proof.

\begin{proposition}\label{ziglin}
Assume system (\ref{ae}) has $k(k\geq 1)$  functionally independent meromorphic first integrals in a neighborhood of $\Gamma$. Then

(i) The normal variational equations (\ref{anve})
have $k(k\geq 1)$  functionally independent first integrals which are rational functions in $\zeta$.

(ii) Each element of Lie algebra $\mathcal{G}$, as linear vector fields, has $k(k\geq 1)$  functionally independent rational first integrals.

\end{proposition}

Let us give some remarks on Proposition \ref{ziglin}. Observing the dimension of (\ref{anve}) is $n-1$, some readers may  doubt the correctness of Lemma \ref{ziglin} in the case of $k=n-1$.
Indeed, for an $n$-dimensional general differential system $\dot x=f(x,t)$, the maximal number of
\emph{autonomous} functionally independent first integrals is $n-1$, whereas that of \emph{time-dependent}  functionally independent first integrals is $n$, see Chapter 10.5-10.7 in \cite{arnold}.
Eq. (30) is an $(n-1)$-dimensional  differential system and may have $(n-1)$ functionally independent first integrals. To further illustrate Proposition \ref{ziglin},
we give two examples with $k=n-1$.

\begin{example} [$n=2,~k=1$] Consider a two-dimensional system
\begin{align}\label{e3}
\dot x=xy,~\dot y=-y^2-x+1,
\end{align}
which has a first integral $\Phi(x,y)=-3x^2+2x^3+3x^2y^2$ and a non-equilibrium solution $\phi(t)=(0,\tanh t)$.
The variational equations along $\phi(t)$ read
\begin{equation}\label{lnve}
\left( {\begin{array}{*{20}{c}}
{\dot \xi }\\
{ \dot \eta }
\end{array}} \right) = \left( {\begin{array}{*{20}{c}}
\tanh t&0 \\
-1&\tanh t
\end{array}} \right)\left( {\begin{array}{*{20}{c}}
\xi \\
\eta
\end{array}} \right),
\end{equation}
and consequently, the normal variational equations read
$$
\dot \xi=a(t)\xi=\tanh t~ \xi,
$$
which has a first integral $\widetilde\Phi(t,\xi)=(\tanh t^2-1)\xi^2.$
\end{example}

\begin{example}[$n=3,~k=2$] Consider a three-dimensional system
\begin{align}\label{e3}
\dot x=-x+yz,~\dot y=-y-xz,~\dot z=-z+xy
\end{align}
which has two first integrals
$
\Phi_1=(x^2+y^2)/(x^2-z^2),~\Phi_2=(y^2+z^2)/(x^2-z^2)
$
 and a non-equilibrium solution $\phi(t)=(0,0,e^{-t})$.
The normal variational equations along $\phi(t)$ are given by
\begin{equation}\label{nve}
\left( {\begin{array}{*{20}{c}}
{\dot \eta_1 }\\
{ \dot \eta_2 }
\end{array}} \right) = \left( {\begin{array}{*{20}{c}}
-1 &e^{-t} \\
-e^{-t} &-1
\end{array}} \right)\left( {\begin{array}{*{20}{c}}
\eta_1 \\
\eta_2
\end{array}} \right).
\end{equation}
One can check that system (\ref{nve}) has two functionally independent first integrals
\begin{align*}
&\widetilde\Phi_1(\eta_1,\eta_2,t)=e^t\sin(e^{-t})\eta_1-e^t\cos(e^{-t})\eta_2,\\
&\widetilde\Phi_2(\eta_1,\eta_2,t)=e^t\cos(e^{-t})\eta_1+e^t\sin(e^{-t})\eta_2.
\end{align*}
\end{example}

In what follows, we aim to establish results analogous to Proposition \ref{ziglin} for which the choice of invariant tensors is changed from first integrals to $n$-forms(or Jacobian multipliers).

\begin{lemma}\label{no2ve}
Assume system (\ref{ae}) has a meromorphic Jacobian multiplier $J(x)$ in a neighborhood of $\Gamma$. Then the variational equations (\ref{ave})
have a time-dependent Jacobian multiplier $J_{VE}(t,\xi)$ which is a rational function with respect to $\xi$.
\end{lemma}
{\bf Proof:}
Let $\xi(t)=x-\psi(t)$, which satisfies
\begin{align}\label{tran1}
\frac{d\xi}{dt}=\widehat F(t,\xi):=F(\xi+\psi(t))-F(\psi(t)).
\end{align}
By Lemma \ref{gjaco2}, system (\ref{tran1}) has a Jacobian multiplier $\widehat J(t,\xi)=J(\xi+\psi(t))$, which is
meromorphic  with respect to $\xi$. By definitions, we have
\begin{align}\label{10-1}
\partial_t \widehat J+\widehat J\nabla_\xi\cdot \widehat F+\langle \partial_\xi \widehat J,\widehat F\rangle\equiv0.
\end{align}
Observing $\widehat J(t,\xi)=P(t,\xi)/Q(t,\xi)$ for certain functions $P,~Q$ which is holomorphic with respect to $\xi$,
we can rewrite (\ref{10-1}) into
\begin{align}\label{10-2}
Q\partial_tP-P\partial_tQ+\langle Q\partial_\xi P-P\partial_\xi Q,\widehat F\rangle+PQ\nabla_\xi\cdot \widehat F\equiv0.
\end{align}
We expand functions $P,Q$ in the neighborhood of $\Gamma$
\begin{align}\label{10-3}
P=P_m(t,\xi)+O(\|\xi\|^{m+1}),Q=Q_k(t,\xi)+O(\|\xi\|^{k+1}),~~P_m\neq0,~Q_k\neq0,
\end{align}
where $P_m,Q_k$ are  the leading terms of $P,~Q$, that is, the lowest order nonvanishing term
of expansions,  $P_m(or~Q_k)$ is homogenous polynomials of degree $m$(or~$k$) with respect to $\xi$.
Similarly, we expand the vector field $\widehat F$ as
\begin{align}\label{10-4}
\widehat F=A(t)\xi+O(\|\xi\|^{2}),
\end{align}
in the neighborhood of $\Gamma$, where $A(t)$ is defined in (\ref{ave}).
Substituting (\ref{10-3})-(\ref{10-4}) into (\ref{10-2}) and comparing the terms of the lowest order $m+k$, we get
\begin{align}\label{10-5}
Q_k\partial_tP_m-P_m\partial_tQ_k+\langle Q_k\partial_\xi P_m-P_m\partial_\xi Q_k,A(t)\xi\rangle+P_mQ_ktr(A(t)) \equiv0,
\end{align}
or equivalently,
\begin{align}\label{10-6}
\partial_t \left(\frac{P_m}{Q_k}\right)+\langle \partial_\xi \left(\frac{P_m}{Q_k}\right),A(t)\xi\rangle+\left(\frac{P_m}{Q_k}\right)tr(A(t))\equiv0,
\end{align}
where $tr$ is the trace of a matrix, i.e., the sum of entries on the diagonal.
It follows from (\ref{10-6}) that the variational equations (\ref{ave}) admit a Jacobian multiplier $J_{VE}(t,\xi)=P_m/Q_k$, which is rational function with respect to $\xi$.\qed

\begin{lemma}\label{ve2nve}
Assume the variational equations (\ref{ave})
have a time-dependent Jacobian multiplier $J_{VE}(t,\xi)$ which is a rational function with respect to $\xi$. Then
the normal variational equations (\ref{anve})
have a time-dependent Jacobian multiplier $J_{NVE}(t,\zeta)$ which is a rational function with respect to $\zeta$.
\end{lemma}
{\bf Proof:} Due to Lemma \ref{gjaco2}, system (\ref{e7}) has a time-dependent Jacobian multiplier
 $$
 \widehat J_{VE}(t,\eta):={J_{VE}(t,P(t)\eta)}{\det P(t)},
 $$
which is  a rational function with respect to $\eta=(\zeta,\eta_1)$. Hence, by definitions we have
\begin{align}\label{11-1}
\partial_t \widehat J_{VE} +\partial_{\eta_1}\widehat J_{VE}  \langle \alpha(t), \zeta \rangle+\langle \partial_\zeta \widehat J_{VE},C(t)\zeta \rangle+ \widehat J_{VE}tr(C(t))\equiv0,
\end{align}
where we use $tr(B(t))=tr(C(t))$. From the proof of Lemma \ref{no2ve}, we see that
$\widehat J_{VE}(t,\eta)=\widehat P_m(t,\eta)/\widehat Q_k(t,\eta)$ for certain functions $\widehat P_m$ and $\widehat Q_k$ which are polynomials of degree $m,~k$  with respect to $\eta$. Hence, (\ref{11-1}) becomes
\begin{align}\label{11-2}
\widehat Q_k\partial_t\widehat P_m-\widehat P_m\partial_t\widehat Q_k+(\widehat Q_k\partial_{\eta_1}\widehat P_m-\widehat P_m\partial_{\eta_1}\widehat Q_k)\langle \alpha(t), \zeta \rangle+\notag\\
\langle \widehat Q_k\partial_\zeta \widehat P_m-\widehat Q_k\partial_\zeta \widehat P_m,C(t)\zeta\rangle+\widehat P_k\widehat Q_mtr(C(t))\equiv0.
\end{align}
We write
\begin{align}\label{11-3}
\widehat P_m=\sum_{j=0}^{m_1}P_{m,j}(t,\zeta)\eta_1^j,~~m_1\leq m,~~P_{m,m_1}\neq0,
\end{align}
with $P_{m,j}$ being polynomials of degree $m-j$ with respect to $\zeta$.
Similarly, we also write
\begin{align}\label{11-4}
\widehat Q_k=\sum_{j=0}^{k_1}Q_{k,j}(t,\zeta)\eta_1^j,~~k_1\leq k,~~Q_{k,k_1}\neq0,
\end{align}
with $Q_{k,j}$ being polynomials of degree $k-j$ with respect to $\zeta$.
Substituting (\ref{11-3})-(\ref{11-4}) into (\ref{11-2}) and comparing the terms of the highest order $m_1+k_1$ with respect to $\eta_1$, we get
\begin{align}\label{11-5}
\partial_t \left(\frac{P_{m,m_1}}{Q_{k,k_1}}\right)+\langle \partial_\zeta \left(\frac{P_{m,m_1}}{Q_{k,k_1}}\right),C(t)\zeta\rangle+\left(\frac{P_{m,m_1}}{Q_{k,k_1}}\right)tr(C(t))\equiv0,
\end{align}
that is,
the normal variational equations  (\ref{anve}) admit a Jacobian multiplier $J_{NVE}(t,\zeta)=P_{m,m_1}/Q_{k,k_1}$, which is a rational function with respect to $\zeta$.\qed

\begin{lemma}\label{nve2galo}
Assume the normal variational equations (\ref{anve})
have a time-dependent Jacobian multiplier $J_{NVE}(t,\zeta)$ which is a rational function with respect to $\zeta$.
Then each element of the Lie  algebra $\mathcal{G}$ of  the identity component $G^0$ of the differential Galois group $G$ for (\ref{anve}), as linear vector fields, has a common rational Jacobian multiplier.
\end{lemma}
{\bf Proof:} Fixed $t_0$, for any $x\in\mathbb{C}^{n-1}$ we consider the solution
$\phi_t(t_0,x_0)$ of (\ref{anve}). Thanks to Lemma \ref{gjaco1}, we have
\begin{align}\label{12-1}
J_{NVE}(t_0,x_0)=J_{NVE}(t,\phi_t(t_0,x_0))\det \partial_{x_0} \phi_t(t_0,x_0),~~\forall t\geq t_0.
\end{align}
Let $\Phi(t)$  be the fundamental-solution matrix of the normal variational equations (\ref{anve}) satisfying $ \det\Phi(t_0)=Id$. Then we have $\phi_t(t_0,x_0)=\Phi(t)x_0$ and (\ref{12-1}) becomes
\begin{align}\label{12-2}
J_{NVE}(t_0,x_0)=J_{NVE}(t,\Phi(t)x_0)\det \Phi(t),~~\forall t\geq t_0.
\end{align}

For any $\sigma\in G$, we have its  representation $M_\sigma$ through $\sigma(\Phi(t))=\Phi(t)M_\sigma$. Furthermore, taking the group action $\sigma$ on the both sides of (\ref{12-2}) yields
\begin{align*}
J_{NVE}(t_0,x_0)&=\sigma\left(J_{NVE}(t,\Phi(t)x_0)\det \Phi(t)\right)=\sigma(J_{NVE}(t,\Phi(t)x_0))\cdot \sigma(\det(\Phi(t))\\
                 &=J_{NVE}(t,\sigma(\Phi(t))x_0)\cdot \det(\sigma(\Phi(t)),\\
                 &=J_{NVE}(t,\Phi(t)M_\sigma x_0)\cdot \det(\Phi(t)M_\sigma),\\
                  &=J_{NVE}(t_0,M_\sigma x_0)\cdot \det(M_\sigma).
\end{align*}
Denote by $J_{gal}(x)=J_{NVE}(t_0,x)$. Based on above discussions, we see that
the rational function $J_{gal}(x)$ satisfies
\begin{align}\label{12-3}
J_{gal}(x)=J_{gal}(M_\sigma x)\det(M_\sigma), \forall x\in\mathbb{C}^{n-1},~\forall M_\sigma\in G.
\end{align}
On the other hand, recall that arbitrary element $Y\in
\mathcal{G}$ can be viewed as a linear vector field: $x\rightarrow
Y(x):=Yx$, for $x\in \mathbb{C}^{n-1}$, and
$e^{Y t}\in G^0$ for all $t\in \mathbb{C}$. It follows from (\ref{12-3}) that
$J_{gal}(x)=J_{gal}(e^{Yt} x)\det (e^{Yt})$ for any $Y\in \mathcal{G}$, which means $J_{gal}(x)$ is a common Jacobian multiplier of a family of linear vector fields $\{\dot x=Yx,~Y\in\mathcal{G}\}$.\qed

Now we can state our main results.
\begin{theorem}\label{thbb}
Let ${\bf \psi}(t)$ be a  non-equilibrium analytic solution  of an $n$-dimensional  analytical system (\ref{ae}) of differential equations. If system (\ref{ae}) has $0\leq k<n$ meromorphic first integrals $\Phi_1(x),\cdots,\Phi_k(x)$ and $n-1-k$
meromorphic Jacobian multipliers $J_1(x),\cdots,J_{n-1-k}(x)$ such that
$$
\Phi_1,\cdots,\Phi_k,\frac{J_2}{J_1},\cdots,\frac{J_{n-1-k}}{J_1}
$$
are functionally independent in a neighborhood of $\Gamma$. Then the following statements hold.

{\rm (i)} The identity
component of the differential Galois group of the  normal variational equations along $\psi(t)$ is abelian.

{\rm (ii)} The identity
component of the differential Galois group of the  variational equations along $\psi(t)$
is solvable.
\end{theorem}
{\bf Proof}: We first show that $J_i/J_1$, $i=2,\cdots,n-1-k$ are first integrals of system (\ref{ae}).
Indeed, by definitions we have
\begin{align}
&J_i(x)\nabla\cdot F(x)+\langle \partial_xJ_i,F\rangle\equiv0,\label{2-1}\\
&J_1(x)\nabla\cdot F(x)+\langle \partial_xJ_1,F\rangle\equiv0\label{2-2}.
\end{align}
Eliminating the divergence term $\nabla\cdot F(x)$ from (\ref{2-1})-(\ref{2-2}) yields
$$
J_1\langle \partial_xJ_i,F\rangle-J_i\langle \partial_xJ_1,F\rangle\equiv0,
$$
which can be rewritten as
$$
\langle \partial_x\left(\frac{J_i}{J_1}\right),F\rangle\equiv0.
$$
Furthermore, $J_i/J_1$ can not be constants, otherwise $
\Phi_1,\cdots,\Phi_k,{J_2}/{J_1},\cdots,{J_{n-1-k}}/{J_1}
$ will be functionally independent. Hence, $J_i/J_1$, $i=2,\cdots,n-1-k$ are first integrals of system (\ref{ae}).

Since (\ref{ae}) has $n-2$  meromorphic first integrals $
\Phi_1,\cdots,\Phi_k,{J_2}/{J_1},\cdots,{J_{n-1-k}}/{J_1}$, by Proposition \ref{ziglin} the
differential Galois group $G$ of the normal variational equations (\ref{anve})
 has $n-2$ rational invariants, and then its Lie algebra $\mathcal{G}$
 has $n-2$ rational first integrals, denoted by $I_1, \cdots, I_{n-2}$.
    Let $U\subset \mathbb{C}^{n-1}$ be a neighborhood of $0$, such that
$I_1, \cdots, I_{n-2}$ are functionally independent on it.
Set the level surface
$$
H_c=\{\eta| I_i(\eta)=c_i, \eta\in \mathbb{C}^{n-1},\  i=1, \cdots,
m\},
$$
where $c_i\in \mathbb{C}$  are constants such that $H_c$ is
regular and $H_c\cap U\neq\emptyset$. Obviously, $H_c$ is an
 one dimensional manifold, and for any fixed point $x\in H_c$, the tangent space of $H_c$ corresponding to $x$  is  one dimensional linear space containing   $\mathcal{V}_{x}=\{Y(x)|Y\in\mathcal{G}\}$.
Therefore, there exists a non-trivial element $Y_0\in \mathcal{G}$ such that
for any linear vector field $Y\in \mathcal{G}$, there exists a one-to-one rational function $a(x)$ such that
$Y=a(x)Y_0$.

Similarly,  since (\ref{ae}) has a meromorphic Jacobian multiplier $J_1$,
by Lemma \ref{nve2galo}, elements $Y(x)$ of the Lie algebra $\mathcal{G}$ have common rational Jacobian multiplier $J_{gal}(x)$, i.e.,
\begin{align}
&J_{gal}\nabla\cdot Y+\langle \partial_xJ_{gal},Y\rangle=J_{gal}(
\langle\partial_x a,Y_0\rangle+a\nabla\cdot Y_0)+a\langle \partial_xJ_{gal},Y_0\rangle\equiv0.\label{2-20}
\end{align}
In particular, for $a(x)=1$ we also have
\begin{align}
J_{gal}\nabla\cdot Y_0+\langle \partial_xJ_{gal},Y_0\rangle\equiv0.\label{2-3}
\end{align}
Combing (\ref{2-20}) and (\ref{2-3}) yields $Y_0(a):=\langle\partial_x a,Y_0\rangle=0$,
that is, $a(x)$ is a rational first integral of the linear vector field $Y_0$.

Consider any two elements $Y_1=a(x)Y_0,~Y_2=b(x)Y_0$ of $\mathcal{G}$ with $a({x}),b({x})\in\mathbb{C}({x})$.
We have $Y_0(a)=Y_0(b)=0$. Then their Lie bracket
$$
[Y_1,Y_2]=[aY_0,bY_0]=(aY_0(b)-bY_0(a))Y_0-ab[Y_0,Y_0]=0,
$$
which means the Lie algebra $\mathcal{G}$ is abelian. Hence the Lie group $G^0$ is also abelian.

Finally, an abelian group is also a
solvable group, and the identity component $\widehat G^0$ of the  variational equations (\ref{e7}) is an extension of the identity component of the normal variational equations $G^0$ (\ref{anve}) by an algebraic group isomorphic to some additive group. This means $G^0$ is solvable if and only if $
\widehat G^0$ is solvable. We complete the proof.\qed

\begin{corollary}
If system (\ref{ae}) has $n-1$ functionally independent meromorphic Jacobian multipliers $J_1(x),\cdots,J_{n-1}(x)$, then the identity
component of the differential Galois group of the  normal variational equations along $\psi(t)$ is abelian,
and that of variational equations along $\psi(t)$ is solvable.
\end{corollary}
{\bf Proof}: It follows from Theorem \ref{thbb} and the fact that the  functional independence of $J_i,i=1,\cdots,{n-1}$ implies the functional independence of $J_i/J_1,i=2,\cdots,{n-1}$.\qed

\begin{corollary}\label{coro2}
Assume system (\ref{ae}) is divergence-free, i.e., $div(F)=0$.  If system (\ref{ae}) has $n-2$ functionally independent meromorphic first integrals $\Phi_1(x),\cdots,\Phi_{n-2}(x)$, then the identity
component of the differential Galois group of the  normal variational equations along $\psi(t)$ is abelian, and that of variational equations along $\psi(t)$ is solvable.
\end{corollary}
{\bf Proof}: Recall that  $div(F)=0$ means system (\ref{ae}) has a Jacobian multiplier $J(x)=1$.\qed

\begin{remark}
In practical applications, as pointed out in Section 4.2 of \cite{ramis} or Section 5.2 of \cite{morales2001},
if the  variational equations are  non-Fuchsian, that is, it has  irregular singularities, one should
  extend $\Gamma$ to a new bigger Riemann surface $\overline{\Gamma}$ by adding possible equilibrium points and  points
at infinity, and should treat the variational equations in a \emph{small} field $\mathcal{M}{\overline{\Gamma}}$ such as $\mathbb{C}(t)$. Then meromorphic non-integrability near  $\overline{\Gamma}$ gives rise to rational non-integrability near $\Gamma$.

\end{remark}

At the end of this section,  we provide two examples to further illustrate Theorem \ref{thbb}.

\begin{example}
Consider the integrable stretch-twist-fold flow
\begin{align}\label{st}
\dot{x}=-8 x y,~~
\dot{y}=11 x^{2}+3 y^{2}+z^{2}+\beta x z-3,~~\dot{z}=2 y z-\beta x y,
\end{align}
which is proposed to model the stretch-twist-fold mechanism of the magnetic field generation
and  plays an important role in understanding the fast dynamo action for magnetohydrodynamics
\cite{st2}.
It admits a polynomial first integral $\Phi=x^3(x^2+y^2+z^2-1)^4$ and
a particular solution $(x(t),y(t),z(t))=(0,-\tanh(3t),0)$. Then
the  variational equations of system (\ref{st}) along this solution are given by
\begin{equation}\label{stve}
\left( {\begin{array}{*{20}{c}}
{\dot \xi }\\
{\dot \eta }\\
{\dot \zeta }
\end{array}} \right) = \left( {\begin{array}{*{20}{c}}
{ 8\tanh(3t)}&0 &0\\
0&{-6\tanh(3t)}&0\\
\beta \tanh(3t)&0&{ -2\tanh(3t) }
\end{array}} \right)\left( {\begin{array}{*{20}{c}}
\xi \\
\eta \\
\zeta
\end{array}} \right).
\end{equation}
According to
Theorem \ref{thbb}, the identity
component of the differential Galois group of the  variational equations (\ref{stve})  along this solution
is solvable, that is, (\ref{stve}) can be integrable by quadrature. Indeed, one can easily
check that
$$
x(t)=c_1\cosh^{8/3}(3t),~~y(t)=\frac{c_2}{\cosh(6t)+1},
~~z(t)=\frac{c_3}{\cosh^{2/3}(3t)}+\frac{c_1\beta}{10}\cosh^{8/3}(3t),
$$
 are the general solutions of  (\ref{stve}).
\end{example}

\begin{example}
Theorem \ref{thbb} can be applied to study the non-existence of first integrals for some high-dimensional divergence-free systems. As a simple
example,
consider  a dynamical system with ${\bf V_4}$ symmetry group \cite{v4}
\begin{align}\label{v4}
\dot{x}=x-yz,~~
\dot{y}=-2y+xz,~~\dot{z}=z-xy,
\end{align}
which has a particular solution $(x(t),y(t),z(t))=(0,0, e^{t})$. The  variational equations read
\begin{equation}\label{v4ve}
\left( {\begin{array}{*{20}{c}}
{\dot \xi }\\
{\dot \eta }\\
{\dot \zeta }
\end{array}} \right) = \left( {\begin{array}{*{20}{c}}
1&-e^{t} &0\\
e^{t}&{-2}&0\\
0&0& 1
\end{array}} \right)\left( {\begin{array}{*{20}{c}}
\xi \\
\eta \\
\zeta
\end{array}} \right).
\end{equation}
The equations for $(\xi,\eta)$ consist of a subsystem
\begin{equation}\label{llnve}
\left( {\begin{array}{*{20}{c}}
{\dot \xi }\\
{\dot \eta }
\end{array}} \right) = \left( {\begin{array}{*{20}{c}}
{ 1}&-e^{t} \\
{e^t}&{ -2}
\end{array}} \right)\left( {\begin{array}{*{20}{c}}
\xi \\
\eta
\end{array}} \right).
\end{equation}
We make a time variable change $\tau=e^t$ to transform (\ref{llnve}) into
\begin{equation}\label{lnve}
\frac{{\rm d}}{{\rm d}\tau}\left( {\begin{array}{*{20}{c}}
{\xi }\\
{ \eta }
\end{array}} \right) = \left( {\begin{array}{*{20}{c}}
{ \frac{1}{\tau}}&-1 \\
{1}&{ -\frac{2}{\tau}}
\end{array}} \right)\left( {\begin{array}{*{20}{c}}
\xi \\
\eta
\end{array}} \right),
\end{equation}
which is equivalent to the Bessel equation
\begin{align}\label{w}
\tau^2\frac{d^2\xi}{d\tau^2}+\tau\frac{d\xi}{d\tau}+(\tau^2-n^2)\xi=0,
\end{align}
with $n=1$. It is well known that the Bessel equation (\ref{w}) has  Liouvillian solutions if and only if $n+1/2$ is an integer, see Chapter 2.8.2 in \cite{ramis}. Hence, we see that the
identity is not solvable.
Suppose system (\ref{v4}) has a rational first integral.
By Corollary  \ref{coro2}, the identity component of the differential Galois group of (\ref{v4ve}) and (\ref{w}) are solvable,
This leads to a contradiction.
Therefore system (\ref{v4}) has no any rational first integrals.
\end{example}

\section{Application to the stationary gravity wave problem in finite depth}

Many studies have been devoted to the phenomenon of stationary gravity waves, which can
 lead to the formation of various
patterns such as solitary waves, star-shaped waves and hexagons waves.
 In particular, Witting \cite{wi} proposed  a new formal series solution of the water wave problem when he studied the solitary wave in a fluid of finite depth. Witting's method has the
advantage of employing a systematic procedure and can give rise to higher approximations for the water waves.
Karabut \cite{ka1,ka2,ka3} showed
the problem of exact summation of this series can be reduced to solve or integrate some homogeneous ordinary differential equations, called Karabut systems.
The aim of this section is to study Karabut systems from the integrability standpoint with the help of Theorem
\ref{thbb}.

\subsection{The mathematical formulation}
To motivate  Karabut
systems to be studied, we give a brief account of the modeling. Consider a plane, incompressible, irrotational stationary flow of a heavy fluid over a flat bottom.
In the coordinate system $(X,Y)$, the origin is located on the bottom, the $X$-axis is directed along the bottom from the left to the right,  the
$Y$-axis
is directed vertically
upward and the free-surface equation is given by $Y=\eta(X)$,  see Figure 1.
Without loss of generality, we assume the fluid flows from the
left to the right. Then we obtain the flow area
$$
\mathcal{D}_1=\left\{Z=X+iY:-\infty<X<+\infty, 0\leq Y\leq \eta(X)\right\}.
 $$
Denote by  $\Phi=\Phi(X,Y)$ the velocity potential and  $\Psi=\Psi(X,Y)$ the streamline function.
Then the flow in the domain $\mathcal{D}$ is  to be potential with the complex potential $\Lambda=\Phi+i\Psi$,
where $i=\sqrt{-1}$ is the  imaginary unit.
In addition, we also denote by $u_0$ and $h_0$ the velocity and the depth
at some point of the free surface. For a solitary wave problem, the point is located at infinity.
For simplicity, we  make a dimensionless rescaling $\chi=\phi+i\psi=\theta(\Phi+i\Psi)/h_0u_0$ in the strip
$$
\mathcal{D}_2=\left\{\chi=\phi+i\psi:-\infty\leq\phi\leq+\infty, 0\leq \psi\leq \theta\right\}.
$$
Here, the Stokes parameter $\theta$ is related to the Froude number $Fr:=u_0/\sqrt{gh_0}$ via $Fr=\sqrt{\tan \theta/\theta}$.
\begin{figure}[htpp]
\centering
\subfigure[Physical flow region .]{
\label{Fig.sub.1}
\includegraphics[width=0.42\textwidth]{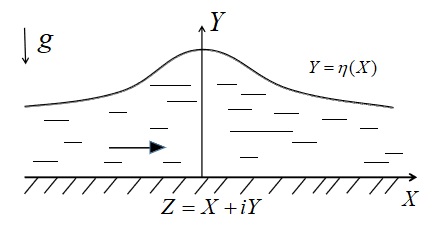}}
\subfigure[Plane of the complex potential $\mathcal{X}$]{
\label{Fig.sub.2}
\includegraphics[width=0.5\textwidth]{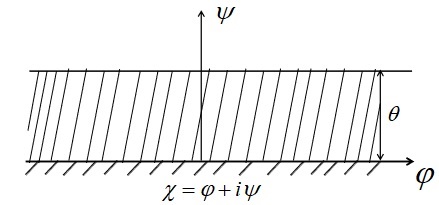}}
\label{Fig.lable}
\end{figure}

Now the  water wave problem is reduced to find  an analytical function $Z=h_0f(\chi)/\theta$  which maps the domain
 $\mathcal{D}_2$ conformally onto the  domain $\mathcal{D}_1$. Set the function $f(\chi)$ in the form $f(\chi)=\chi+W(\chi)$.
 At the free surface $Y=\eta(X)$, by the Bernoulli integral, we have
 \begin{align*}
\frac{1}{2}\left|\frac{\mathrm{d}(\Phi+i \Psi)}{\mathrm{d} Z}\right|^{2}+g Y=\frac{1}{2} u_{0}^{2}+g h_{0},
 \end{align*}
which is equivalent to
\begin{align}\label{w1}
\bigg|\frac{dW}{d\chi}+1\bigg|^2=\frac{1}{1-2\nu {\rm Im}~W},~~\textit{for}~~\psi=\theta
\end{align}
with $\nu=\cot \theta$. At the flat bottom, we have
\begin{align}\label{w2}
 {\rm Im}~W=0,~~\textit{for}~~\psi=0.
\end{align}
Moreover, for the solitary wave, the boundary condition at infinity
\begin{align}\label{w3}
\mathop {\lim }\limits_{\varphi \to -\infty } {\mathop{\rm Im}\nolimits}~W = \mathop {\lim }\limits_{\varphi \to +\infty } {\mathop{\rm Im}\nolimits}~W = 0
\end{align}
should be fulfilled. In a word, the mathematical statement of the bottom, sides and upper boundary
condition is given in (\ref{w1})-(\ref{w3}).

The shallow-water theory shows that it is
reasonable to consider $W(\phi+i\psi)$ as a periodic function of variable $\psi$.
For small-amplitude waves, one can consider the
solution of (\ref{w1})-(\ref{w3}) with the form
\begin{align}\label{ddw}
W=\sum_{j=1}^{\infty}{\theta^{2j}W^{(j)}(\chi)}
\end{align}
 with  $W^{(j)}(\chi)$ being polynomials of
$\cosh(\chi/2)^{-2}$, which corresponds to the shallow-water expansion.
Witting \cite{wi} proposed a new expansion parameter for solitary waves and constructed a solution of (\ref{w1})-(\ref{w2}) with the form
\begin{align}\label{ddw2}
W=\sum_{j=1}^{\infty}{E_j(\theta)\zeta^j},~~{\rm Im}~E_j=0,~~\zeta=e^{\chi},
\end{align}
which permits  calculation to extremely high order. When $\theta=\pi/n$ with $n$ being an odd integer, Karabut \cite{ka1,ka2,ka3} introduced the new unknown functions
$$
P_j(\chi)=W(\zeta e^{i(2j-2)\theta}),~~~~j=1,2,\cdots,n
$$
and proved that $P_j$ satisfy the following system of ordinary differential equations
\begin{align}\label{p}
 \bigg(\frac{dP_{j+1}}{d\chi}+1\bigg) \bigg(\frac{dP_{j}}{d\chi}+1\bigg)=\frac{1}{f_j},~~~~j=1,2,\cdots,n,
\end{align}
where $P_{n+1}=P_1,~~f_j=1+i\nu(P_{j+1}-P_{j}).$ Clearly, to solve the boundary-value problem (\ref{w1})-(\ref{w3})
in the form of the Witting  series, one only needs to  integrate system (\ref{p}) and to take $W = P_1$.

Rewrite equations (\ref{p}) into equations  with respect to the variables $f_j$
\begin{align}\label{pf}
 \frac{d f_j}{d\chi}=\mathrm{i}\nu\frac{\prod\limits_{k=1}^{(n-1)/2}{f_{[(2k+j-1)\mod n]+1}}-\prod\limits_{k=1}^{(n-1)/2}{f_{[(2k+j-2)\mod n]+1}}}{\sqrt{\prod\limits_{k = 1}^n {{f_k}} }}, ~~j=1,2,\cdots,n.
\end{align}
Making a time scale $dt=i\nu/\sqrt{\prod\limits_{k = 1}^n {{f_k}} }d\chi$ and replacing $f_j$ with $x_j$, (\ref{pf}) become the following equivalent homogeneous equations
\begin{align}\label{pf2}
 \frac{d x_j }{dt}=\prod_{k=1}^{(n-1)/2}{x_{[(2k+j-1)\mod n]+1}}-\prod_{k=1}^{(n-1)/2}{x_{[(2k+j-2)\mod n]+1}},~~~~j=1,2,\cdots,n,
\end{align}
which are called Karabut systems.

\subsection{Integrability analysis of the 3-dimensional Karabut systems}

For $\theta=\pi/3$, we obtain the 3-dimensional Karabut system
\begin{align}\label{k3}
\left\{
\begin{aligned}
\dot x_1&=x_3-x_2,\\
\dot x_2&=x_1-x_3,\\
\dot x_3&=x_2-x_1.
\end{aligned}
\right.
\end{align}
Obviously, system (\ref{k3}) is a linear system of differential equations and one can get its general   solutions,
which implies the Witting series for $\theta=\pi/3$ can be summed up exactly \cite{ka1}. From the point of view of integrability, we have the following results.
\begin{proposition}\label{pro1}
System $(\ref{k3})$ is completely integrable with two functionally independent first integrals
$$
I_1=x_1x_2+x_2x_3+x_3x_1,~~~~~~~I_2=x_1+x_2+x_3.
$$
\end{proposition}

Using first integrals $I_1$ and $I_2$, we can construct infinitely many Hamilton-Poisson realizations
$(\mathbb{R}^3, \Pi_{a,b}, H_{c,d})$ of system (\ref{k3})
parameterized by the group $SL(2,\mathbb{C})$, which implies system (\ref{k3}) is bi-Hamiltonian.
\begin{proposition}\label{pro2}
Let $a,b,c,d\in \mathbb{R}$ such that $ad-bc=1$.
For system (\ref{k3}), there exist infinitely many  Hamilton-Poisson realizations
$(\mathcal{P},\Pi_{ab},H_{cd})$ such that
$$
\dot{\bf x}=\Pi_{a,b}\cdot \nabla H_{c,d},~~{\bf x}=(x_1,x_2,x_3)^T,
$$
 where the Poisson structure $\Pi_{a,b}$ and the Hamiltonian $H_{c,d}$ are given by
\begin{align*}
\Pi_{a,b}=&\left[\begin{array}{*{20}{c}}
0&{ - b - ax_1 - ax_2}&{b + ax_1 + ax_3}\\
{b + ax_1 + ax_2}&0&{ - b - ax_2 - ax_3}\\
{ - b - ax_1 - ax_3}&{b+ ax_2 + ax_3}&0
\end{array} \right],\\[1mm]
H_{c,d}=&c(x_1x_2+x_2x_3+x_3x_1)+d(x_1+x_2+x_3),
\end{align*}
here $a, b, c, d$ satisfy
$\left[ {\begin{array}{*{20}{c}}
a&b\\
c&d
\end{array}} \right]
\in SL(2,\mathbb{R})$.
\end{proposition}

Next result shows that the 3-dimensional Karabut system admits a Lax formulation.
\begin{proposition}\label{pro3}
The system (\ref{k3}) can be written in the Lax form $\dot L= [N, L]$, where the matrices $L$ and respectively $N$ are given by
$$
L:=\left(
     \begin{array}{ccc}
       0 & -x_3 & x_2 \\
       x_3 & 0 & -x_1 \\
        -x_2& x_1 & 0 \\
     \end{array}
   \right),~~~
   N:=\left(
     \begin{array}{ccc}
       0 & -1-x_3 & 1+x_2 \\
       1+x_3& 0 & -1-x_1 \\
        -1-x_2& 1+x_1 & 0 \\
     \end{array}
   \right).
$$
\end{proposition}

\subsection{Integrability analysis of the 5-dimensional Karabut systems}

For $\theta=\pi/5$,  we obtain the 5-dimensional Karabut system
\begin{align}\label{k5}
\left\{
\begin{aligned}
\dot x_1&=x_3x_5-x_2x_4,\\
\dot x_2&=x_4x_1-x_3x_5,\\
\dot x_3&=x_5x_2-x_4x_1,\\
\dot x_4&=x_1x_3-x_5x_2,\\
\dot x_5&=x_2x_4-x_1x_3.
\end{aligned}
\right.
\end{align}
Karabut \cite{ka3} found two functionally independent polynomial first integrals of (\ref{k5}),
$$
\Phi_1=x_1+x_2+x_3+x_4+x_5,~~~~\Phi_2=x_1^2+x_2^2+x_3^2+x_4^2+x_5^2.
$$
But as he said,
\emph{ it is not known whether the system has other polynomial
integrals}, so that he had to study system (\ref{k5}) by numerical methods.  The goal of this section is to
deal with this problem.  More specifically, we
 investigate its  polynomial integrability, namely, \emph{what is the maximal
number of functionally independent polynomial first integrals that system $(\ref{k5})$
can exhibit?}

A nice result is due to Christov \cite{k1}, in which he used a Morales-Ramis type theory on Bogoyavlenskij integrability \cite{zung} and showed that
system (\ref{k5}) is not meromorphically integrable in the Bogoyavlenskij sense.
Since we don't know if (\ref{k5}) admits other symmetry fields besides the trivial  symmetry field $w_1=F$,  we can only conclude from Christov's result that the number of functionally independent polynomial first integrals
of system (\ref{k5}) is less than four. Otherwise, (\ref{k5}) has four independent polynomial first integrals, so it is
integrable in the Bogoyavlenskij sense, which yields a contradiction.
For this corollary, there is a  simple proof: observing the Karabut system is a  quadratic homogeneous polynomial differential system and we can get its Kovalevskaya exponents $(-1, 1,2,3/2+\sqrt{-7-8i}/2, 3/2-\sqrt{-7-8i}/2)$ corresponding to a balance $c=(-i,-1,1,i,0)$, where $i=\sqrt{-1}$. Then by Theorem 5.5 in \cite{gg},
 system (\ref{k5}) has no four functionally independent polynomial first integrals. The non-rational  Kovalevskaya exponents also imply
that  system (\ref{k5}) does not enjoy the Painlev\'{e}/weak-Painlev\'{e}  property. In addition, we mention that the analysis of differential Galois groups in Christov's work \cite{k1} is obtained by computing  local monodromy matrices. As our aim is to give a practical example to show the effectiveness of Theorem \ref{thbb}, we will provide a systematic procedure to analyze the differential Galois group associated with Karabut systems with the help of Kovacic's algorithm, which can be applied to other high-dimensional systems.

\begin{theorem}\label{mth}
System $(\ref{k5})$ has two and only two functionally independent meromorphic  first integrals.
\end{theorem}

A direct consequence is that system (\ref{k5}) has no additional polynomial first integrals which is
functionally independent with $\Phi_1$ and $\Phi_2$.

\noindent{\bf Proof of Theorem \ref{mth}:}
The basic idea of our proof is to find an
integrable invariant manifold $\mathcal{N}$ for the 5-dimensional Karabut system (\ref{k5}), to reduce the (normal) variational equations  along a solution contained in $\mathcal{N}$ into a second-order equation with rational coefficients, and to prove that the identity component of the differential
Galois group associated with this second-order equation is not solvable.

System (\ref{k5}) admits an \emph{integrable invariant manifold} defined by
$$
\mathcal{N}:=\big\{(x_1, x_2, x_3, x_4, x_5)\in\mathbb{C}^5:x_1+x_4=0,~x_2+x_3=0,~x_5=0\big\}.
$$
Indeed, system (\ref{k5}) restricted to $\mathcal{N}$ is given by
\begin{align}\label{e2}
\dot x_1=x_1x_2,~~~\dot x_2=-x_1^2,
\end{align}
which is completely integrable with a polynomial first integral $\Phi=x_1^2+x_2^2$.
Given the initial condition $(x_1(t_0),x_2(t_0))=(x_{10},x_{20})\neq(0,0)$,  we get
$$
x_2=\pm\sqrt{x_{10}^2+x_{20}^2-x_1^2},~~\pm\frac{dx_1}{x_1\sqrt{x_{10}^2+x_{20}^2-x_1^2}}=dt,
$$
that is,
$$
x_1(t)=\frac{4\exp(\frac{t+C_2}{C_1})}{4+C_1^2\exp(\frac{2t+2C_2}{C_1})},~~x_2(t)=\frac{\dot x_1}{x_1},
$$
or
$$
x_1(t)=\frac{4\exp(\frac{t+C_2}{C_1})}{C_1^2+4\exp(\frac{2t+2C_2}{C_1})},~~x_2(t)=\frac{\dot x_1}{x_1},
$$
where $C_1=(x_{10}^2+x_{20}^2)^{-1/2}$ and $C_2$ is an integration constant.
For convenience, we set $C_1=1,C_2=0$ and get a non-equilibrium solution of system (\ref{k5})
$$
\phi(t)=\left(\frac{4\exp(t)}{\exp(2t)+4},\frac{4-\exp(2t)}{\exp(2t)+4},\frac{\exp(2t)-4}{\exp(2t)+4},-\frac{4\exp(t)}{\exp(2t)+4},0\right).
$$

Let $\Gamma$ be the phase curve corresponding to the solution $\phi(t)$. Then the variational equations along $\Gamma$ read
\begin{align}\label{ve}
\frac{{\rm d}}{{\rm d}t}{\bf \xi}={\bf A}(t){\bf \xi},
\end{align}
where
$$
{\bf A}(t)=\left( {\begin{array}{*{20}{c}}
{0}&{\frac{4\exp(t)}{\exp(2t)+4}}&{0}&{\frac{\exp(2t)-4}{\exp(2t)+4}}&{\frac{\exp(2t)-4}{\exp(2t)+4}}\\
{-\frac{4\exp(t)}{\exp(2t)+4}}&{0}&{0}&{\frac{4\exp(t)}{\exp(2t)+4}}&{\frac{4-\exp(2t)}{\exp(2t)+4}}\\
{\frac{4\exp(t)}{\exp(2t)+4}}&{0}&{0}&{-\frac{4\exp(t)}{\exp(2t)+4}}&{\frac{4-\exp(2t)}{\exp(2t)+4}}\\
{\frac{\exp(2t)-4}{\exp(2t)+4}}&{0}&{\frac{4\exp(t)}{\exp(2t)+4}}&{0}&{\frac{\exp(2t)-4}{\exp(2t)+4}}\\
{\frac{4-\exp(2t)}{\exp(2t)+4}}&{-\frac{4\exp(t)}{\exp(2t)+4}}&{-\frac{4\exp(t)}{\exp(2t)+4}}&{\frac{4-\exp(2t)}{\exp(2t)+4}}&{0}
\end{array}} \right).
$$
Since system (\ref{k5}) has two first integrals $\Phi_1$ and $\Phi_2$, by Lemma \ref{ss}, system (\ref{ve}) has two time-dependent first
integrals
\begin{align}\label{fir}
G_1(t,\xi)=\xi_1+\xi_2+\xi_3+\xi_4+\xi_5,~G_2(t,\xi)=\frac{8\exp(t)}{\exp(2t)+4}(\xi_1-\xi_4)+\frac{8-2\exp(2t)}{\exp(2t)+4}(\xi_2-\xi_3).
\end{align}
Then, by means of the change of independent variable
\begin{align}\label{tr}
t\rightarrow \tau:=\exp(t)-4\exp(-t),
\end{align}
we transform (\ref{ve}) into the linear differential system with rational coefficients
\begin{align}\label{ve2}
\frac{{\rm d}}{{\rm d}\tau}{\bf \xi}={\bf B}(\tau){\bf \xi},
\end{align}
where
$$
{\bf B}(t)=\left( {\begin{array}{*{20}{c}}
{0}&{\frac{4}{16+\tau^2}}&{0}&{\frac{\tau}{16+\tau^2}}&{\frac{\tau}{16+\tau^2}}\\
{-\frac{4}{16+\tau^2}}&{0}&{0}&{\frac{4}{16+\tau^2}}&{-\frac{\tau}{16+\tau^2}}\\
{\frac{4}{16+\tau^2}}&{0}&{0}&{-\frac{4}{16+\tau^2}}&{-\frac{\tau}{16+\tau^2}}\\
{\frac{\tau}{16+\tau^2}}&{0}&{\frac{4}{16+\tau^2}}&{0}&{\frac{\tau}{16+\tau^2}}\\
{-\frac{\tau}{16+\tau^2}}&{-\frac{4}{16+\tau^2}}&{-\frac{4}{16+\tau^2}}&{-\frac{\tau}{16+\tau^2}}&{0}
\end{array}} \right).
$$

It should be pointed out that the transformation (\ref{tr}) does not change the
identity component of the differential galois group (see  Theorem 2.5 in \cite{ramis}).
Making a linear transformation $\xi=P\eta$ with
\begin{align}\label{pp}
P=\left( {\begin{array}{*{20}{c}}
{1}&{0}&{0}&{0}&{-1}\\
{0}&{1}&{0}&{-1}&{0}\\
{0}&{0}&{0}&{1}&{0}\\
{0}&{0}&{0}&{0}&{1}\\
{0}&{0}&{1}&{0}&{0}
\end{array}} \right),
\end{align}
system (\ref{ve2}) becomes the following equivalent form
\begin{align}\label{ve3}
\frac{{\rm d}}{{\rm d}\tau}{\bf \eta}={\bf C}(t){\bf \eta},
\end{align}
where
$$
{\bf C}(t)=\left( {\begin{array}{*{20}{c}}
{\frac{\tau}{16+\tau^2}}&{\frac{4}{16+\tau^2}}&{\frac{2\tau}{16+\tau^2}}&{0}&{0}\\
{0}&{0}&{-\frac{2\tau}{16+\tau^2}}&{0}&{0}\\
{-\frac{\tau}{16+\tau^2}}&{-\frac{4}{16+\tau^2}}&{0}&{0}&{0}\\
{\frac{4}{16+\tau^2}}&{0}&{-\frac{\tau}{16+\tau^2}}&{0}&{-\frac{8}{16+\tau^2}}\\
{\frac{\tau}{16+\tau^2}}&{0}&{\frac{\tau}{16+\tau^2}}&{\frac{4}{16+\tau^2}}&{-\frac{\tau}{16+\tau^2}}
\end{array}} \right).
$$
Then we obtain a $3$-dimensional subsystem

\begin{align}\label{3e}
 \frac{{\rm d}}{{\rm d}\tau}\left( {\begin{array}{*{20}{c}}
{\eta_1}\\
{\eta_2}\\
{\eta_3}
\end{array}} \right) = \left( {\begin{array}{*{20}{c}}
{\frac{\tau}{16+\tau^2}}&{\frac{4}{16+\tau^2}}&{\frac{2\tau}{16+\tau^2}}\\
{0}&{0}&{-\frac{2\tau}{16+\tau^2}}\\
{-\frac{\tau}{16+\tau^2}}&{-\frac{4}{16+\tau^2}}&{0}\\
\end{array}} \right)\left( {\begin{array}{*{20}{c}}
{\eta_1}\\
{\eta_2}\\
{\eta_3}
\end{array}} \right).
\end{align}
Due to (\ref{fir}) and (\ref{pp}), one can easily check that $H(\eta_1,\eta_2,\eta_3)=\eta_1(t)+\eta_2(t)+\eta_3(t)$ is a first integral of (\ref{3e}), that is,
$$
\frac{{\rm d}}{{\rm d}t}\left(\eta_1(t)+\eta_2(t)+\eta_3(t)\right)=0,
$$
or $\eta_1(t)+\eta_2(t)+\eta_3(t)=$ constant. Let this constant be zero, then by (\ref{3e}), we obtain
\begin{align}\label{2e}
\frac{{\rm d}}{{\rm d}\tau}
\left( {\begin{array}{*{20}{c}}
{\eta_1}\\
{\eta_2}
\end{array}} \right) = \left( {\begin{array}{*{20}{c}}
{-\frac{\tau}{16+\tau^2}}&{\frac{4-2\tau}{16+\tau^2}}\\
{\frac{2\tau}{16+\tau^2}}&{\frac{2\tau}{16+\tau^2}}
\end{array}} \right)\left( {\begin{array}{*{20}{c}}
{\eta_1}\\
{\eta_2}
\end{array}} \right).
\end{align}
From the first equation of (\ref{2e}), we have
\begin{align}\label{eta2}
\eta_2=\frac{16+\tau^2}{4-2\tau}\frac{{\rm d}\eta_1}{{\rm d}\tau}+\frac{\tau\eta_1}{4-2\tau}.
\end{align}
Substituting (\ref{eta2}) into the second equation of (\ref{2e}), we can eliminate the variable $\eta_2$ and
get an equivalent second-order equation
\begin{align}\label{2ee}
\frac{{\rm d}^2}{{\rm d}\tau^2}\eta_1+P(\tau)\frac{{\rm d}}{{\rm d}\tau}\eta_1+Q(\tau)\eta_1=0,
\end{align}
where the coefficients of (\ref{2ee})
are as follows
\begin{align*}
P=-\frac{2\tau+16}{(\tau-2)(\tau^2+16)},~~Q=\frac{2\tau^3-14\tau^2+16\tau-32}{(\tau+16)^2(\tau-2)}.
\end{align*}
Next, under the change of dependent variable
$$
\eta_1=\chi \exp\left[-\frac{1}{2}\int P{\rm d}\tau\right],
$$
(\ref{2ee}) is converted to its reduced form
\begin{align}\label{ren}
\frac{{\rm d}^2}{{\rm d}\tau^2}\chi=r(\tau)\chi,
\end{align}
with
\begin{align*}
r=\frac{3}{4(\tau-2)^2}-\frac{11-8i}{16(\tau+4i)^2}-\frac{11+8i}{16(\tau-4i)^2}
+\frac{1}{20(\tau-2)}-\frac{8+59i}{320(\tau+4i)}-\frac{8-59i}{320(\tau-4i)}.
\end{align*}

We claim that the identity component of the differential Galois group of system (\ref{ren}) is not
solvable. To this end, we use the Kovacic's algorithm \cite{k}. Obviously, (\ref{ren})
is Fuchsian with four regular singular points $\tau=2,4i,-4i,\infty$, all of them are of order two.
For case 1 in Kovacic's algorithm, by simple computations, we get
$$
\alpha_\infty^\pm=\frac{1}{2}\pm\frac{\sqrt{7}i}{2},~~~~~\alpha_2^+=\frac{3}{2},~~~~~~~\alpha_2^-=\frac{1}{2},
$$
$$
\alpha_{4i}^\pm=\frac{1}{2}\pm\frac{\sqrt{-7-8i}}{4},~~~~~~~~~~~~~
\alpha_{-4i}^\pm=\frac{1}{2}\pm\frac{\sqrt{-7+8i}}{4},
$$
then $d:=\alpha^\pm_\infty-\alpha^\pm_2-\alpha^\pm_{4i}-\alpha^\pm_{-4i}\in\mathbb{C}/\mathbb{R}$ is not a non-negative integer, which means
case 1 of Lemma \ref{kth} is impossible. For case 2, we can obtain the auxiliary sets $E_i$ from the Kovacic's Algorithm
$$
E_\infty=E_{4i}=E_{-4i}=\{2\},~~~~E_2=\{-2,2,6\}.
$$
Next, we should select elements $e_\infty\in E_{\infty}$, $e_{\pm4i}\in E_{\pm4i}$ and $e_2\in E_{2}$ such that
$d:=(e_\infty-e_{4i}-e_{-4i}-e_2)/2$ is a non-negative integer. There
is only one possible choice of $(e_\infty,e_{4i},e_{-4i},e_2)=(2,2,2,-2)$ with $d=0$. We construct the rational function
$$
\theta(\tau)=\frac{-1}{\tau-2}+\frac{2\tau}{\tau^2+16},
$$
and check if there exists a monic polynomial $P$ of zero degree satisfying the equation
\begin{align}\label{1q}
P'''+3\theta
P''+(3\theta^2+3\theta'-4r)P'+(\theta''+3\theta\theta'+\theta^3-4r\theta-2r')P=0.
\end{align}
However,
$$
\theta''+3\theta\theta'+\theta^3-4r\theta-2r'=\frac{8\tau^2+52\tau+224}{(\tau-2)^2(\tau^2+16)^2}\neq0.
$$
Hence, case 2 of Lemma \ref{kth} is also excluded. For case 3, by Theorem 2.1 in \cite{k}, it is also impossible that the differential Galois group is finite since
the necessary conditions cannot hold. In summary, the differential Galois group of (\ref{ren}) falls into case 4  of Lemma \ref{kth} and is
$SL(2,\mathbb{C})$. In this case, the identity component is also $SL(2,\mathbb{C})$, which is not solvable.

Assume system $(\ref{k5})$ has three functionally independent meromorphic  first integrals.
Due to Corollary \ref{coro2}, the identity component of the differential Galois group of (\ref{ve}),
(\ref{ve2}) or its equivalent form (\ref{ve3}) is solvable. Thanks to Lemma \ref{sub}, so the subsystem (\ref{3e}) is, i.e., system (\ref{3e}) can be solved by quadrature.
So system (\ref{2ee}) is also solvable. Finally, note that (\ref{2ee}) and (\ref{ren}) have the same Liouvillian solvability, which implies  (\ref{ren}) is solvable, namely, the identity component of (\ref{ren}) is solvable. This leads to a contradiction.  \qed

\begin{remark}
 Similar to $3D$ and $5D$ Karabut systems,  Karabut systems with dimension $n\geq 7$ also have two first integrals
 $\Phi_1=x_1+\cdots x_n$ and $\Phi_2=x_1^2+\cdots x_n^2$, but their integrability analysis is still open.
The main difficulty is to find an analytic  particular solution lying in an integrable invariant
manifold.  The numerical simulations in \cite{ka3} inspire us to believe that Karabut systems with dimension $n\geq 7$ are also non-integrability.
\end{remark}



\appendix
\renewcommand{\appendixname}{Appendix}

\section*{Acknowledgements}
S.Shi was partially supported  by the National Natural
Science Foundation of China (No. 11771177), China Automobile Industry Innovation
and Development Joint Fund (No. U1664257), Program for Changbaishan Scholars of Jilin Province
and Program for JLU Science, Technology Innovative Research Team (No. 2017TD-20). K. Huang was partially supported by
Sichuan University postdoctoral interdisciplinary Innovation Fund (No.0020104153010),
the Fundamental Research Funds for the Central Universities (No.20826041E4168),
the National Natural
Science Foundation of China (No.12001386, No.12090013).


\begin{thebibliography}{00}


 \bibitem{ga11}
Acosta-Hum\'{a}nez P B, \'{A}lvarez-Ram\'{\i}rez M, Delgado J.
Non-integrability of some few body problems in two degrees of freedom.
{Qual Theory Dyn Syst}, 2009, 8: 209-239.







 \bibitem{ga7}
Acosta-Hum\'{a}nez P B, \'{A}lvarez-Ram\'{\i}rez M, Stuchi T J.
Nonintegrability of the Armbruster-Guckenheimer-Kim quartic Hamiltonian through Morales-Ramis theory.
{SIAM J Appl Dyn Syst}, 2018, 17: 78-96.


\bibitem{poly}
Acosta-Hum\'{a}nez P B, Bl\'{a}zquez-Sanz D.
Non-integrability of some Hamiltonians with rational potentials.
{Discrete Contin Dyn Syst Ser B}, 2008, 10: 265-293.


 \bibitem{zung}
Ayoul M, Zung N T. Galoisian obstructions to non-Hamiltonian integrability.
C R Math Acad Sci Paris, 2010, 348: 1323-1326


\bibitem{arnold}
Arnold V I. Ordinary Differential Equations. Berlin: Springer-Verlag, 1992.

\bibitem{st2}
Bajer K, Moffatt H K. On a class of steady confined Stokes flows with chaotic streamlines.
{J Fluid Mech}, 1990, 212: 337-363.





\bibitem{non1}
Bountis T. Investigating non-integrability and chaos in complex time. {Phys D}, 1995, 86: 256-267.


\bibitem{non2}
Bolsinov A V, Taimanov I A. Integrable geodesic flows with positive topological entropy. {Invent Math},
2000, 140: 639-650.


\bibitem{B}
Bogoyavlenskij O I. Extended integrability and bi-hamiltonian systems. {Comm Math Phys}, 1998, 196: 19-51.



\bibitem{baider1996}
Baider A, Churchill R C, Rod D L, Singer M F. On the infinitesimal geometry
of integrable systems. Mech Day, 1996, 7: 5-56.


\bibitem{churchill}
Churchill R C, Rod D L, Singer M F. Group-theoretic obstructions to integrability.
{Ergod Theory  Dyn Syst}, 1995, 15: 15-48.


 \bibitem{jacobi1}
Casale G. Morales-Ramis theorems via Malgrange pseudogroup.
{Ann Inst Fourier (Grenoble)}, 2009, 59: 2593-2610.



\bibitem{k1}
Christov O. Non-integrability of the Karabut system.
{Nonlinear Anal Real World Appl}, 2016, 32: 91-97.



\bibitem{jordan}
Duval G, Maciejewski A J.
Jordan obstruction to the integrability of Hamiltonian systems with homogeneous potentials.
{Ann Inst Fourier (Grenoble)}, 2009, 59: 2839-2890.



\bibitem{gg}
Goriely A.
Integrability and nonintegrability of dynamical systems. Singapore: World Scientific, 2001.



\bibitem{hu2014}
Hu Y X. On the first integrals of n-th order autonomous systems.
{J Math Anal Appl}, 2018, 459: 1062-1078.



\bibitem{in1}
Jiao J, Huang K Y, Liu W S, Stationary Shear Flows of Nematic Liquid Crystals: A Comprehensive Study via Ericksen-Leslie Model. {J Dyn Diff Equat}, 2022, 34: 239-269.



\bibitem{jo0}
Jacobi C G J. Vorlesungen \"{u}ber Dynamik. Berlin: Druck und Verlag von G. Reimer, 1884.

\bibitem{jo1}
Jovanovi\'{c} B. Geometry and integrability of Euler-Poincar\'{e}-Suslov equations.
{Nonlinearity}, 2001, 14: 1555-1567.




\bibitem{jo2}
Kozlov V V. Symmetries, Topology and Resonances in Hamiltonian Mechanics. Berlin: Springer-Verlag, 1996.





 \bibitem{kozlov2013}
Kozlov V V. The Euler-Jacobi-Lie integrability theorem. {Regul Chaotic Dyn}, 2013, 18: 329-343.



\bibitem{ka1}
Karabut E. Summation of the Witting series in a problem on a solitary wave. {Siberian Math J}, 1995, 36: 287-304.

\bibitem{ka2}
Karabut E. Asymptotic expansions in the solitary-wave problem. J Fluid Mech, 1996, 318: 109-123.

\bibitem{ka3}
Karabut E. On the summation of Witting series in the solitary wave problem. {J Appl Mech Tech Phys}, 1999, 40: 36-45.

\bibitem{k}
Kovacic J J.
An algorithm for solving second order linear homogeneous differential equations.
J Symbolic Comput, 1986, 2: 3-43.

\bibitem{js2016}
Llibre J, Valls C, Zhang X. The completely integrable differential systems are essentially linear
differential systems. {J Nonlinear Sci}, 2015, 25: 815-826.


\bibitem{in2}
Llibre J, Tian Y. Dynamics of the FitzHugh-Nagumo system having invariant algebraic surfaces.
 {Z Angew Math Phys}, 2021, 72: 15, 17pp.


\bibitem{in3}
Llibre J, Messias M.
Global dynamics of the Rikitake system.
{Phys D}, 2009, 238: 241-252.



 \bibitem{li2012}
Li W L, Shi S S. Galoisian obstruction to the integrability of general dynamical systems.
{J Differential Equations}, 2012, 252: 5518-5534.


\bibitem{v4}
Letellier C, Gilmore R.
Symmetry groups for 3D dynamical systems.
{J Phys A}, 2007, 40: 5597-5620.


\bibitem{ga4}
Maciejewski A J, Przybylska M.
All meromorphically integrable 2D Hamiltonian systems with homogeneous potential of degree 3.
Phys Lett A, 2004, 327: 461-473.

 \bibitem{li2008}
Maciejewski A J, Przybylska M, Yoshida Y. Necessary conditions for super-integrability of Hamiltonian systems.
{Phys Lett A}, 2008, 372: 5581-5587.




 \bibitem{jacobi2}
Maciejewski A J, Przybylska M. Integrability analysis of the stretch-twist-fold flow.
{J Nonlinear Sci}, 2020, 30: 1607-1649.





\bibitem{ramis}
Morales-Ruiz J J. Differential Galois theory and non-integrability of Hamiltonian
systems. Basel: Birkh\"{a}user Verlag, 1999.


\bibitem{morales2001}
Morales-Ruiz J J, Ramis J P. Galoisian obstructions to integrability of Hamiltonian
systems I.  Methods~Appl Anal, 2001, 8: 33-95.


\bibitem{morales1994picard}
Morales-Ruiz J, Sim{\'o} C. Picard-vessiot theory and Ziglin's
  theorem. J Differential Equations, 1994, 107: 140-162.

 \bibitem{ga2}
Morales-Ruiz J J, Sim{\'o} C, Simon S. Algebraic proof of the non-integrability of
Hill's problem. {Ergodic Teory Dynam Systems}, 2005, 25: 1237-1256.

 \bibitem{ga3}

Maciejewski A J, Przybylska M.
Darboux points and integrability of Hamiltonian systems with homogeneous polynomial potential.
{J Math Phys}, 2005, 46: 062901, 33pp.


\bibitem{ga5}
Maciejewski A J, Przybylska M, Stachowiak T. Nonexistence of the final first integral in the Zipoy-Voorhees space-time. {Physical Review D}, 2013, 88: 064003.

 \bibitem{ga6}
Maciejewski A J, Przybylska M, Yaremko Y.
Dynamics of a dipole in a stationary electromagnetic field.
{Proc A}, 2019, 475: 20190230, 20 pp.



 \bibitem{jacobi3}
Przybylska M. Differential Galois obstructions for integrability of homogeneous
Newton equations. {J Math Phys}, 2008, 49: 022701, 40pp.



 \bibitem{ga71}
Szumi\'{n}ski W, Wo\'{z}niak D. Dynamics and integrability analysis of two pendulums coupled by a spring.
{Commun Nonlinear Sci Numer Simul}, 2020, 83: 105099, 16 pp.


\bibitem{singer1993}
Singer M F, Ulmer F.
Galois groups of second and third order linear differential equations.
{J Symbolic Comput}, 1993, 16: 9-36.




 \bibitem{ga1}
Tsygvintsev A. The meromorphic non-integrability of the three-body problem.
{J Reine Angew Math}, 2001, 537: 127-149.


\bibitem{van}
Van der Put M, Singer M F. Galois Theory of Linear
Differential Equations. Berlin Heidelberg: Springer-Verlag, 2003.




\bibitem{wi}
Witting J. On the highest and other solitary waves. {J Appl Math}, 1975, 28: 700-719.


\bibitem{non3}
Y. Yagasaki, Galoisian obstructions to integrability and Melnikov criteria for chaos in two-degree-of-freedom
Hamiltonian systems with saddle centres. {Nonlinearity}, 2003, 16: 2003-2012.




\bibitem{zx}
Zhang X. Integrability of Dynamical Systems: Algebra and Analysis. Singapore: Springer-Verlag, 2017.



 \bibitem{zhang2014}
Zhang X. Liouvillian integrability of polynomial differential systems. {Trans Amer Math Soc}, 2016, 368: 607-620.



\bibitem{ziglin}
Ziglin S L. Branching of solutions and nonexistence of first integrals in Hamiltonian mechanics I.
 {Funct Anal Appl}, 1983, 16: 181-189.


















\end{thebibliography}
\end{document}